\documentclass[12pt]{article}
\usepackage{amsmath,amsfonts}
\usepackage{amsthm}
\usepackage{amssymb}
\usepackage{amscd}
\usepackage{xypic}
\usepackage{verbatim}
\usepackage[plainpages=false,colorlinks,hyperindex,pdfpagemode=None,bookmarksopen,linkcolor=red,citecolor=blue,urlcolor=blue]{hyperref}
\usepackage{pdflscape}
\usepackage{stmaryrd}
\usepackage{mathabx}
\usepackage{multirow}
\usepackage{appendix}
 \swapnumbers
\theoremstyle{plain}

\renewcommand{\marginpar}[1] {  }
\renewcommand{\comment}[1] {  }

\usepackage[all]{xy}
\catcode`\é=\active\def é{\'e}
\catcode`\è=\active\def è{\`e}
\catcode`\à=\active\def à{\`a}
\catcode`\ù=\active\def ù{\`u}

\catcode`\ê=\active\def ê{\^{e}}
\catcode`\î=\active\def î{\^{i}}
\catcode`\ô=\active\def ô{\^{o}}
\catcode`\û=\active\def û{\^{u}}
\catcode`\ç=\active\def ç{\c{c}}

\newtheorem{theo}{Theorem}[section]
\newtheorem{lem}[theo]{Lemma}

\newtheorem{cor}[theo]{Corollary}

\theoremstyle{remark}

\numberwithin{equation}{section}

\def \g{\mathfrak}

\def\bb{\backslash}

\def\al{\alpha}
\def\be{\beta}

\def\De{\Delta}
\def\de{\delta}
\def\ep{\varepsilon}

\def\ga{\gamma}

\def\La{\Lambda}
\def\la{\lambda}
\def\om{\omega}
\def\Om{\Omega}

\def\si{\sigma}
\def\Si{\Sigma}
\def\Te{\Theta}
\def\te{\theta}

\def\bb{\backslash}
\def\ad{{\rm{ad}}}
\def\Ad{{\rm{Ad}}}
\def\det{{\rm{det}}}

\def\cD{{\mathcal D}}

\def\cG{{\mathcal G}}
\def\cH{{\mathcal H}}

\def\cL{{\mathcal L}}
\def\cM{{\mathcal M}}

\def\cP{{\mathcal P}}

\def\cS{{\mathcal S}}
\def\cT{{\mathcal T}}

\def\cY{{\mathcal Y}}
\def\cZ{{\mathcal Z}}
\DeclareFontFamily{OT1}{rsfs}{}
\DeclareFontShape{OT1}{rsfs}{n}{it}{<-> rsfs10}{}
\DeclareMathAlphabet{\mathscr}{OT1}{rsfs}{n}{it}
\def\Lr{{\mathscr L}}

\def\E{{\rm E}}
\def\F{{\rm F}}

\def\C{\mathbb C}

\def\R{\mathbb R}

\def\N{\mathbb N}

\def\R{\mathbb R}
\def\Z{\mathbb Z}

\def\scP{\underline{\mathcal P}}
\def\scZ{\underline{\mathcal Z}}

\def\sA{\underline{ A}}

\def\sG{\underline{ G}}
\def\sH{\underline{ H}}

\def\sJ{\underline{ J}}

\def\sM{\underline{ M}}
\def\sN{\underline{ N}}

\def\sP{\underline{ P}}
\def\sp{\underline{ p}}

\def\sS{\underline{ S}}
\def\sT{\underline{ T}}

\def\sp{\underline{ p}}
\def\bi{\bigskip}
\def\me{\medskip}
\def\no{\noindent}

\def\ste{\par\smallskip\noindent}
\def\ste{\par\smallskip\noindent}

\def\dem{ {\em Proof~: \ste }}
 \def\beq{\begin{equation}}
\def\eeq{\end{equation}}

\newenvironment{res}
               {\begin{equation}\begin{minipage}{0.85\textwidth}}
               {\end{minipage}\end{equation}}
\def\ber{\begin{res}}
\def\eer{\end{res}}

\def\qed{{\null\hfill\ \raise3pt\hbox{\framebox[0.1in]{}}\break\null}}

\begin{document}

\author{P. Delorme\thanks{The first author was supported by a grant of Agence Nationale de la Recherche with reference
ANR-13-BS01-0012 FERPLAY.}, P. Harinck, S. Souaifi}
\title{Geometric side of a  local  relative trace formula}
\date{}
\maketitle
\maketitle

\begin{abstract} Following a scheme suggested by B. Feigon, we investigate a  local relative trace formula in the situation of a reductive $p$-adic
group $G$ relative to a symmetric subgroup $H= \underline{H}(F)$ where $\underline{H}$ is split over the local field $\F$ of characteristic zero
 and $G = \underline{G} (F)$ is the restriction of scalars of  $\underline{H} _{I \E}$ relative to a quadratic
unramified extension $\E$ of $\F$. We  adapt techniques of the proof of the local trace formula by  J.~Arthur  in  order to get a  geometric expansion of the integral  over $H \times H$ of a  truncated kernel associated to the regular representation of $G$.
\end{abstract}

 \medskip


\tableofcontents
\addcontentsline{toc}{section}{Introduction} \section*{Introduction}  In this article, we investigate a local relative trace formula in the situation of $p$-adic groups relative to a symmetric subgroup. This work is inspired by the recent results of B. Feigon 
(\cite{F}), where she investigated what she called a local relative trace formula on $PGL(2)$ and a local Kuznetsov trace formula for $U(2)$.

Before we describe our setting and results, we would to explain on the toy model of finite groups the framework of the formulas of B.~Feigon.
We even start with  the more general framework of the relative trace formula initiated by H. Jacquet (\cite{Ja}).

Let $G$ be a finite group and let $H$, $H'$, $\Gamma$ be subgroups of $G$.
We endow any finite set with the counting measure. We denote by $r$ the   right regular representation of $G$ on $L^2(\Gamma \bb G)$ and we consider the
 $H$-fixed linear form $\xi$  on  $L^2(\Gamma \bb G)$ defined by
\beq \label{xi} \xi= \sum _{h \in H\cap \Gamma \bb H} \delta_{\Gamma h} \eeq
where $\delta_{\Gamma h}$ is the Dirac measure of the coset $\Gamma h$, 
or in other words $$\xi(\psi)= \int_{H\cap \Gamma \bb H } \psi(\Gamma h) dh,\quad \psi\in L^2(\Gamma\bb G).$$
We define similarly $\xi'$ relative to  $H'$.\\
We view $\xi$, $\xi'$ as elements of $L^2(\Gamma \bb G)$ and we form the coefficient $c_{\xi, \xi'}(g)= ( r(g)\xi,  \xi')$. Integrating over functions on $G$, it defines a ''distribution'' $\Theta$ on $G$ which is right invariant  by $H$ and left invariant by $H'$.  
The relative trace formula in this context gives two expressions of $\Theta(f) $ for $f$ a function on $G$, the first one, called the geometric side, in terms of orbital integrals, and the second one, called the spectral side, in terms of irreducible representations of $G$.\medskip

First we deal with the geometric side. For this purpose we introduce   suitable orbital integrals. For $\gamma \in \Gamma$, we set $[\gamma]:= (H'\cap \Gamma)  \gamma  (H \cap  \Gamma )$ and one introduces two subgroups of $H'\times H$ $$(H'\times H)_{\gamma}=\{(h',h)\vert h'\gamma h^{-1}= \gamma \}, (H'\cap \Gamma\times H\cap\Gamma)_{\gamma}= (H'\times H)_{\gamma} \cap (\Gamma\times \Gamma ).$$
Then, we define the orbital integral of a function $f$ on $G$ by
 $$I([\gamma], f)= \int_{(H'\times H)_{\gamma} \bb (H'\times H)} f(h'\gamma h^{-1})dh'dh.$$
Let $f$ be a function on $G$. Since $r(g)\delta_{\Gamma h}=\delta_{\Gamma hg^{-1}}$, the definition of $\xi$ and $\xi'$ gives

$$ \Theta(f)=\sum_{g\in G} f(g)\Theta(g) =\sum_{g\in G} f(g)\frac{1}{vol(\Gamma\cap H)} \frac{1}{vol(\Gamma\cap H')}\sum_{h\in H} \sum_{h'\in H'}(\delta_{\Gamma hg^{-1}}, \delta_{\Gamma h'}).$$

Changing  $g$ in $g^{-1}h$ and using the fact that  $( \delta_{\Gamma g},  \delta_{\Gamma h'})$ is equal to $1$  for $g\in \Gamma h'$ and to zero otherwise, one gets
\beq\label{Teta1}\Theta (f)= \frac{1}{vol(\Gamma\cap H)} \frac{1}{vol(\Gamma\cap H')}\sum_{h\in H} \sum_{h'\in H'}\sum_{\gamma\in \Gamma}  f(h'\gamma h).\eeq

A simple computation of volumes leads to the geometric expression of $\Theta$ in terms of orbital integrals
 \beq \label{Tgeo} \Theta(f)=  \sum_{[\gamma ]\in H'\cap \Gamma  \bb \Gamma  / \Gamma \cap H}  vol ( (H'\cap \Gamma \times H \cap \Gamma)_\gamma\bb (H'\times H)_\gamma) I([\gamma], f).\eeq 
 
Let us turn to the spectral side. We decompose $L^2(\Gamma\bb G)$ into isotypic components $\oplus _{\pi\in \hat{G}}\cH_\pi$
The restriction of $\xi$ and $\xi'$ to $\cH_\pi$ will be denoted $\xi_\pi$ and $\xi'_\pi$ respectively. The spectral formula for $\Theta$ is the simple equality
\beq\label{Tspec}\Theta= \sum_{\pi\in \hat{G}} c_{\xi_\pi, \xi'_\pi}.\eeq
Notice that it might be  also interesting to decompose further the representation into irreducible representations and the restriction of $\xi$ to each of them will be called a period.\medskip

There is a third interpretation of the distribution $\Theta$. If $f$ is a function on $G$, then the operator $r(f)$ on $L^2(\Gamma \bb G)$ is an integral operator whose kernel $K_f$ is the function on $\Gamma\bb G\times \Gamma\bb G$ given by 
   $$K_f(x,y)=\sum_{\gamma\in\Gamma} f(x^{-1}\gamma y).$$\\
By (\ref{Teta1}), one gets easily the following expression of $\Theta(f)$  
   \beq \label{Kf} \Theta( f)= \int_{ (H'\cap \Gamma \bb H') \times (H\cap \Gamma \bb H} K_{f} (h', h) dh' dh. \eeq
This point of view is probably the best one. But it is important to have the representation theoretic meaning of $\Theta$.\medskip
 
 The toy model for the local relative trace formula of B.~Feigon appears as a particular case of the above relative trace formula. 
In that case, the groups  $G$, $H$ and $H'$  are   products $G_1 \times G_1$, $ H_1\times H_1$ and $ H'_1\times H'_1$ respectively and $\Gamma$ is the diagonal of $G_1\times G_1$. 
Then $\Gamma \bb G$ identifies with $G_1$ and the right representation corresponds to the  representation $R$ of $G_1\times G_1$ on $L^2(G_1)$ given by
$[R(x,y)\phi](g)=\phi(x^{-1}gy).$
Then, we have 
$$\xi (\psi) = \int_{H_1} \psi (h)dh, \quad \psi\in L^2(G_1).$$

The spectral side is more concrete. If $(\pi_1, \cH_{\pi_1})$ is an irreducible unitary  representation of $G_1$ then $G_1\times G_1$ acts on $\mbox{End} (\cH_{\pi_1})$ by an irreducible  representation denoted by $\pi$. It is unitary if we use the scalar product associated to the Hilbert-Schmidt norm.  Moreover $L^2(G_1) $ is canonically isomorphic to the direct sum $\oplus_{\pi_1\in \hat{G_1}} \mbox{End} (\cH_{\pi_1})$. Let $P_\pi$ be the orthogonal projector onto the space of invariant vectors under $H_1$. Then,
the period map $\xi_\pi$, which  is a linear form on $\mbox{End} (\cH_{\pi_1})$, is  given by
$$\xi_\pi (T) = \int_{H_1} Tr( \pi_1(h)  T)dh=(T, P_\pi),\quad T\in \mbox{End} (\cH_{\pi_1}).$$

One further decomposes $\xi_\pi$ by using an orthonormal basis ($\eta_{\pi_1, i}$) of the space of $H_1$-invariant vectors. We will use the identification 
of $\mbox{End}  (\cH_{\pi_1})$ with the tensor product of $\cH_{\pi_1}$ with its conjugate complex vector space. In this identification, one has
$$P_\pi =\sum_i \eta_{\pi_1, i} \otimes \eta_{\pi_1, i}.$$ 
We define similar notations for $\xi'$ relative to $H'$. 
 Then, for two functions $f_1, f_2 $ on $G_1$, the spectral side (\ref{Tspec}) can be written
$$\Theta (f_1 \otimes f_2)=  \sum_{\pi_1 \in \hat{G_1}} \sum_{i ,i'} c_{\eta_{\pi_1, i}, \eta '_{\pi_1, i'}}(f_1)c_{\eta_{\pi_1, i}, \eta'_{\pi_1, i'}}(f_2).$$

For the geometric side,   we define the integral orbital of a function $f$ on $G_1$  by   
 $$I(g, f)= \int_{(H'_1 \times H_1)_g \bb H'_1 \times H_1} f(h'gh^{-1}) dh dh'$$
 which depends only on the double coset $H'_1 g H_1$.
Then one gets by (\ref{Tgeo}) the equality
$$\Theta(f_1\otimes f_2)  = \sum_{g \in H'_1 \bb G_1/ H_1}   v(g ) I( g , f_1) I(g  , f_2) $$
where the $v(g )$'s are positive constants depending  on volumes.
Hence the final form of the local relative trace formula is:
$$\sum_{g \in H'_1 \bb G_1/ H_1}   v(g ) I( g, f_1) I(g, f_2) =  \sum_{\pi_1 \in \hat{G_1}} \sum_{i ,i'} c_{\eta_{\pi_1, i}, \eta '_{\pi_1, i'}}(f_1)c_{\eta_{\pi_1, i}, \eta'_{\pi_1, i'}}(f_2).$$

This formula allows to invert the orbital integrals $I(g, f_1)$. For this purpose, one chooses $g_1\in G_1$ and takes for $f_2$   the Dirac measure at $g_1$. Then $I(g_1, f_2)= 1 $ and the other orbital integrals of $f_2$ are zero.
Hence $$ v(g_1) I(g_1, f_1)=  \sum_{\pi_1 \in \hat{G_1}} \sum_{i ,i'} c_{  \eta_{\pi_1, i}, \eta '_{\pi_1, i'} }(f_1)c_{\eta_{\pi_1, i}, \eta'_{\pi_1, i'}}(f_2). $$
In order to make the formula more precise, one needs to compute the constants $c_{\eta_{\pi_1, i}, \eta'_{\pi_1, i'}} (f_2)$.\medskip

The inversion of orbital integrals is one of our motivations to investigate  a local relative trace formula in the situation of $p$-adic groups relative to a symmetric subgroup $H$ and we will take $H=H'$. \medskip

In this article, we consider a reductive algebraic group $\sH$ defined over  a non archimedean local field $\F$  of characteristic  $0$. We fix   a quadratic unramified extension $\E$ of $\F$ and we consider the group    $\sG:=\mbox{Res}_{\E/\F} \sH$  obtained by restriction of scalars of    $\sH$, where here $\sH$ is  considered as a group defined over $\E$. We denote by $H$ and $G$ the group  of $\F$-points of $\sH$ and $\sG$ respectively. Then  $G$ is isomorphic to $\sH(E)$ and $H$ appears as the fixed points of $G$ under the involution of $G$ induced by the nontrivial element of the Galois group of $\E/\F$.  We assume that $\sH $ is split over $\F$ and we fix a maximal split torus $A_0$ of $H$. The groups $G$ and $H$ correspond to $G_1$ and $H_1=H_1'$ respectively in our example of a	 local relative trace formula for finite groups.\me

The starting point of our study is the analogue to the expression  (\ref{Kf}).  We consider the regular representation $R$ of $G\times G$ on $L^2(G)$ given by $(R(g_1,g_2)\psi)(x)=\psi(g_1^{-1} xg_2)$.  Then for $f=f_1\otimes f_2$ where $f_1$ and $f_2$ are two smooth compactly supported functions on $G$, the corresponding operator $R(f)$ is an integral operator on $L^2(G)$ with smooth kernel
$$K_{f }(x,y) =\int_G f_1(xg) f_2(gy) dg= \int_G f_1(g) f_2(x^{-1} g y) dg.$$

As $H$ may  be not compact, even modulo the split component $A_H$ of the center of $H$,  we have to truncate this kernel to integrate it. We multiply this kernel   by a product of functions $u(x, T) u(y,T)$ where $u(\cdot, T)$ is the characteristic function of a large compact subset in $A_H\bb H$ depending on a parameter $T\in a_0=\mbox{Rat}(A_0)\otimes_\Z\R$ ($\mbox{Rat}(A_0)$ is the group of $\F$-rational characters of $A_0$)  as in \cite{ArLT}  (cf.   (\ref{uxT})). As $H$ is   split, we have $A_H=A_G$. Hence the kernel $K_f$ is left invariant by the diagonal  $diag(A_H)$ of $A_H$ and we can integrate the truncated kernel over $diag(A_H)\bb H\times H$. We set 
$$K^T(f):=\int_{diag(A_H)\bb (H\times H)}  K_{f}(x_1, x_2) u(x_1,T) u(x_2,T) d\overline{(x_1, x_2)}.$$

In \cite{ArLT}, J.~Arthur studies the integral of  $K_f(x,x) u(x,T)$ over $A_G\bb G$ to obtain its local trace formula on reductive groups.\me

We study the geometric expression of the distribution $K^T(f)$ and its dependence on the parameter $T$. Our main results (Theorem \ref{resultatprincipal} and Corollary \ref{expressionJT})  assert that 
 $K^T(f)$ is asymptotic as $T$ approaches infinity to another distribution $J^T(f)$ of the form 
\beq\label{asymp} J^T(f)=\sum_{k=0}^N p_{\xi_k}(T,f) e^{\xi_k(T)}\eeq
where   $\xi_0=0,\ldots \xi_N$ are distinct points of the dual space $i a_0^*$ and  each $p_{\xi_k}(T,f)$ is a polynomial function  in $T$. Moreover, the constant term $\tilde{J}(f):=p_0(0,f)$ of $J^T(f)$ is well-defined and  uniquely determined by $K^T(f)$. We give an explicit expression of this constant term in terms of weighted orbital integrals.

These results are analogous to those of \cite{ArLT} for   the group case. Our proof follows closely the study by J.~Arthur of the geometric side of his local trace formula which we were able to  adapt under our assumptions to the case of double truncations. \me

In the first section, we  introduce notation on groups and  on symmetric spaces according to \cite{RR}.   The starting point of our study is the Weyl integration formula established in \cite{RR}, which takes into account the $(H,H)$-double classes of $\si$-regular elements of $G$ (cf. (\ref{Thm3.4RR}) and (\ref{Weyl1})). These double classes are express in terms of $\si$-torus, that is torus whose elements are anti-invariant by $\si$.  Under our assumptions, there is a bijective correspondence $S\to S_\si$ between maximal tori of $H$ and maximal $\si$-tori of $G$ which preserves $H$-conjugacy classes. 

Then the Weyl integration  formula can be written in terms of Levi subgroups  $M\in\cL(A_0)$ of $H$ containing $A_0$ and $M$-conjugacy classes of maximal anisotropic tori of $M$  (cf. (\ref{Weyl2})):
$$ \int_G f(g) dg=\sum_{M\in\cL(A_0)} c_M\sum_{S\in\cT_M} \sum_{x_m\in \kappa_S}c_{S, x_m}\int_{S_\si} |\De_\si(x_m\ga)|_\F^{1/2} \int_{diag(A_M)\bb H\times H}  f(h^{-1} x_m\ga l)d\overline{(h,l)} d\ga$$
where $\kappa_S$ is a finite subset of $G$, $c_M$ and $ c_{S,x_m}$ are positive constants, $\cT_M$ is a suitable set of anisotropic tori of $M$ and $\De_\si$ is a jacobian. 

A fundamental  result for our proofs concerns  the orbital integral $\cM(f)$ of  a compactly smooth function $f$ on $G$. It is defined  on $\si$-regular points by 
$$\cM(f)(x_m\ga)=|\De_\si(x_m\ga)|_\F^{1/4} \int_{diag(A_S)\bb H\times H}  f(h^{-1} x_m\ga l)d\overline{(h,l)}, $$
where $S$ is a maximal torus of $H$, $x_m\in \kappa_S$ and  $ \ga  \in S_\si $ such that  $x_m\ga$ is $\si$-regular. As in the group case using the exponential map and the property that each root of $S_\si$  has multiciplity $2$ in the Lie algebra of $G$, we prove that the orbital integral 
is bounded on the subset of $\si$-regular  points of $G$ (cf. Theorem \ref{OIbounded}). \me

In the second section, we  explain the truncation process based on the notion of $(H,M)$-orthogonal sets and prove our main results.  Using the Weyl integration formula,  we can write 
$$ K^T(f)=\sum_{M\in\cL(A_0)} c_M\sum_{S\in\cT_M}\sum_{x_m\in\kappa_S} c_{S,x_m}\int_{S_\si} K^T(x_m, \ga, f) d\ga$$
where  $$K^T(x_m, \ga,f)=|\De_\si(x_m\ga)|_\F^{1/2}\int_{diag(A_M)\bb H\times H} \int_{diag(A_M)\bb H\times H} f_1(y_1^{-1}x_m\ga y_2) $$
$$\times f_2(x_1^{-1} x_m\ga   x_2) u_M(x_1,y_1,x_2,y_2, T)   d\overline{(x_1,x_2)} d\overline{(y_1,y_2)}$$ and  
$$u_M(x_1,y_1,x_2,y_2, T) =\int_{A_H\bb A_M} u(y_1^{-1}ax_1,T) u(y_2^{-1}ax_2, T) da.$$
The function $J^T(f)$ is obtained in a similar way to $K^T(f)$ where we replace the weight function $u_M(x_1,y_1,x_2,y_2, T) $ by another weight function $v_M(x_1,y_1,x_2,y_2, T) $.\me

The weight function $v_M$ is given by

$$ v_M(x_1,y_1,x_2,y_2, T) := \int_{A_H\bb A_M}\si_M(h_M(a), \cY_M(x_1, y_1, x_2, y_2, T)) da$$
where $\si_M(\cdot,\cY) $ is 
 the function   of \cite{ArLT} depending on a $(H,M)$-orthogonal set  $\cY$ and $\cY_M(x_1, y_1, x_2, y_2, T)$ is a $(H,M)$- orthogonal set obtained as the "minimum" of two $(H,M)$- orthogonal sets $\cY_M(x_1, y_1, T)$ and $\cY_M(x_2, y_2, T)$ (cf. (\ref{SigmaM}), Lemma \ref{infY} and (\ref{YxyT})). 
 If $\cY_1$ and $\cY_2$ are two $(H,M)$-orthogonal positive sets  then the "minimum" $\cZ$ of $\cY_1$ and $\cY_2$ satisfies the property that the convex hull  $\cS_M(\cZ)$ in $a_H\bb a_M$ of the points of $\cZ$ is the intersection of the convex hulls  $\cS_M(\cY_1)$ and $\cS_M(\cY_2)$ in $a_H\bb a_M$ of the points of $\cY_1$ and $\cY_2$  respectively.
 
  If $\Vert T\Vert$ is  large relative to $\Vert x_i\Vert, \Vert y_i\Vert, i=1,2$ then $\sigma_M(\cdot, \cY_M(x_1, y_1, x_2, y_2, T))$ is just the characteristic function of $\cS_M(\cY_M(x_1, y_1, x_2, y_2, T))$. In that case, this function  is equal to the product of  $\sigma_M(\cdot, \cY_M(x_1, y_1,T))$ and  $\sigma_M(\cdot, \cY_M(x_2, y_2, T))$.
 
 Our proofs consist to establish good estimates of $\vert u_M((x_1,y_1,x_2,y_2, T)-v_M(x_1,y_1,x_2,y_2, T)\vert$  when $x_i,y_i,i=1,2$ satisfy $f_1(y_1^{-1}x_m \ga y_2)f_1(x_1^{-1}x_m \ga x_2)\neq 0$ for some $\ga\in S_\si$ and $x_m\in\kappa_S$. Then, using that orbital integrals are bounded, we deduce our result on $\vert K^T(f)-J^T(f)\vert$. \me
 
 This work is a first step towards a local relative trace formula. For the spectral side, we have to prove that $K^T(f)$ is asymptotic to a distribution $k^T(f)$ which is of general form (\ref{asymp}) and  constructed from spectral data. We hope that we can express the constant term of $k^T(f)$ in terms of regularized local period integrals introduced by B. Feigon in \cite{F} in the same way than Jacquet-Lapid-Rogawski  regularized period integrals for automorphic forms in \cite{JLR}. We plan to explicit such a local relative trace formula for $PGL(2)$.\me
 
 \noindent {\bf Acknowledgments.} We thank warmly Bertrand Lemaire  for his answers to our many questions on algebraic groups. We thank Bertrand Rémy and David Renard for our helpful discussions. We thank also Guy Henniart for providing  us a proof of (\ref{G1}).

\section{ Preliminaries}
\subsection{ Reductive $p$-adic groups}\label{defgroups}

Let  $\F$ be a non archimedean local field of characteristic  $0$ and odd  residual characteristic $q$. Let $|\cdot|_\F$ denote the normalized valuation   on $\F$.

For an   algebraic variety   $\sM$  defined over  $\F$, we identify $\sM$  with   $\sM(\overline{\F})$ where $\overline{\F}$ is an  algebraic closure of $\F$ and we set $M:=\sM(\F).$ 

We will use conventions like in \cite{W2}. 
 One considers various algebraic groups $\sJ$ defined over $\F$,  and   sentences
like 
 \ber \label{corse} " let $M$ be an algebraic  group" will mean '' let $M$ be the $\F$-points of an algebraic group   $\sM$ defined  over $\F$"  and    " let $A$ be a split torus
 "  will mean " let
$A$ be the  group of $\F$-points   of a torus,  $\underline A$,   defined and split  over $\F$ ."\eer 
If  $J$ is an algebraic group, one denotes by $\mbox{Rat} (J)$ 
the  group of its rational characters defined over $\F$. 
If $V$ is a vector space, $V^*$ will denote its dual. If $V$ is real, 
$V_\C$ will denote its complexification.\me

Let $\sG$  be an algebraic reductive group defined over $\F$.
We fix a maximal split torus $A_0$ of $G$ and we denote by $M_0$ its centralizer  in $G$. 

We denote by $A_G$ the maximal split torus of the center of $G$ and we define
$$a_{G}:= {\rm Hom}_{\Z} (\mbox{Rat} (G),\R).$$ \\
One has the canonical map  $h_{G}: G \rightarrow  a_{G}$
 which is defined by 
\beq \label{hG} e^{\langle h_{G}(x), \chi\rangle}= \vert \chi (x)\vert_\F,\quad  x\in G,
\chi
\in \mbox{Rat} (G).\eeq\\
The restriction of rational characters 
 from  $G$ to  $A_{G}$  induces an  isomorphism 
\beq \label{iso}  \mbox{Rat} (G)\otimes _{\Z}\R \simeq  \mbox{Rat} (A_{G})\otimes _{\Z} \R.\eeq

Notice that $\mbox{Rat} (A_G)$ appears as a  generating lattice in the dual space $ a^*_G$ of $a_G$ and  \beq \label{aotimes} a^*_G \simeq \mbox{Rat} (G) \otimes _\Z \R. \eeq   

 The kernel of  
$h_{G}$, which is denoted by $G^1$,   is the intersection of the kernels   of $\vert \chi \vert_{\F}$ for all  character $\chi \in \mbox{Rat}
(G)$  of
$G$. The groupe $G^1$ is distinguished in $G$ and contains  the derived group $G_{der}$ of $G$. Moreover, it is well-known that
\ber\label{G1} the group $G^1$ is generated by  the compact subgroups of $G$.\eer
G. Henniart has communicated to us an unpublished proof of this result   by N. Abe, F. Herzig, G. Henniart and M.–F. Vigneras.

\ber\label{aGF}One denotes by  $a_{G, \F}$ (resp.,  $\tilde{a}_{G, \F}$)  the image of  $G$ (resp.,  $A_{G}$) by  $ h_{G}$. Then 
$G/G^{1}$ is  isomorphic to the lattice 
$a_{G,\F}$.  \eer\ste 

 If $P$ is a parabolic subgroup of $G$ with Levi subgroup $M$,
we  keep the same notation with $M$ instead of $G$.
\me

The inclusions
$ A_{G}\subset  A_{M}\subset M\subset G$  determine a surjective  morphism    $ a_{M, \F}\rightarrow a_{G, \F}$ (resp., an injective  morphism, 
 $ \tilde 
{a}_{G, \F}   \rightarrow \tilde 
{a}_{M, \F}$) which extends uniquely to a surjective linear map  $h_{MG}$ from    $ a_{M}$ to $a_{G}$ (resp., injective map  between
$a_{G}$ and 
$a_{M}$).
The second map allows to identify  $a_{G}$ with a subspace of 
$a_{M}$ and the kernel of the first one, $a^{G}_{M}$,  satisfies 
\ber $$\label{aMG} a_{M}= a^{G}_{M}\oplus a_{G}.$$
\eer
For $M=M_0$, we set $a_0:=a_{M_0}$ and $a_0^G:=a_{M_0}^G$. We fix a scalar product $(\cdot,\cdot)$ on $a_0$ which is invariant under the Weyl group $W(G,A_0)$ of $(G,A_0)$. Then $a_G$ identifies with the fixed point set of $a_0$ by $W(G,A_0)$ and $a_0^G$ is an invariant subspace of $a_0$ under $W(G,A_0)$. Hence, it is the orthogonal subspace to $a_G$ in $a_0$. The space $a_G^*$ might be viewed as a subspace of $a^*_0$ by (\ref{aMG}). Moreover, by definition of the surjective map $a_0\to a_G$, one deduces that
\ber\label{hGM} if $m_0\in M_0$ then  $h_{G}(m_0)$ is the orthogonal projection of $h_{M_0}(m_0)$ onto $a_G$.\eer
From (\ref{aMG}) applied to $(M,M_0)$ instead of $(G,M)$, one obtains a decomposition $a_0=a_0^M\oplus a_M$. From the $W(G,A_0)$ invariance of the scalar product, one gets:
\ber\label{a0M} The decomposition $a_0=a_0^M\oplus a_M$ is an orthogonal decomposition. 

The space $a_M^*$ appears as a subspace of $a_0^*$ and, in the identification of $a_0$ with $a^*_0$ given by the scalar product, $a_M^*$ identifies with $a_M$.\eer

The decomposition $a_M=a_M^G\oplus a_G$ is orthogonal relative to the restriction to $a_M$ of the $W(G,A_0)$-invariant inner product on $a_0$ and the natural map $h_{MG}$ is identified with the orthogonal projection of $a_M$ onto $a_G$. 
\ber\label{projreseau} In particular, $a_{G,\F}$ is the orthogonal projection of $a_{M,\F}$ onto $a_G$. Moreover, we have $\tilde{a}_{G,\F}= a_G\cap \tilde{a}_{M,\F}$ (cf. \cite{ArLT} (1.4)).\eer

By a Levi subgroup of $G$, we mean a group $M$ containing $M_0$ which is the Levi component of a parabolic subgroup of $G$.  If $P$ is a parabolic subgroup containing $M_0$ then it has a unique Levi subgroup  denoted by $M_P$ which contains $M_0$. We will denote by $N_P$ the unipotent radical of $P$.

For a Levi subgroup $M$, we write  $\cL(M)$ for the finite set of  Levi subgroups of $G$ which contain $M$ and we also let $\cP(M)$ denote the finite set of   parabolic subgroups $P$ with $M_P=M$. \me

 Let $K$ be  the fixator of a special point in the apartment of $A_0$ in the Bruhat-Tits building. We have the Cartan decomposition
\beq\label{Cartan}  \quad G=KM_0K.\eeq\\
If $P=M_PN_P$ is a parabolic subgroup of $G$ containing $M_0$, then 
\beq\label{Iwasawa} G= PK=M_PN_PK.\eeq\\
If $x\in G$, we can write 
\beq\label{mP} x= m_P(x)n_P(x)k_P(x),\quad m_P(x)\in M_P, n_P(x)\in N_P, k_P(x)\in K.  \eeq
We set  \beq\label{hP} h_P(x):= h_{M_P}(m_P(x)).\eeq
The point $m_P(x)$ is defined up an element of $K$ but $h_p(x)$ does not depend of this choice.\me

 We  introduce a norm $\Vert\cdot\Vert$ on $G$ as in (\cite{W2} \textsection I.1.) (called height function in (\cite{ArLT})). Let   $\La_0: G\to GL_n(\F)$ be an  algebraic embedding. 
For $g\in G$, we write
$$\La_0(g)=(a_{i,j})_{i,j=1\ldots n},\quad  \La_0(g^{-1})=(b_{i,j})_{i,j=1\ldots n}.$$
We set 
\beq\label{norme}\Vert g\Vert:=\sup_{i,j} \sup(|a_{i,j}|_\F, |b_{i,j}|_\F).\eeq
If $\La: G\to GL_d(\F)$ is another algebraic embedding then  the norm $\Vert\cdot \Vert_\La$ attached to $\La$ as above is equivalent to $\Vert\cdot\Vert$ in the following sense: there are a positive constant $C_\La$ and a positive integer $d_\La$ such that  
$$\Vert g\Vert_\La\leq C_\La \Vert g\Vert^{d_\La}.$$
This allows us to   use results of \cite{W2} for estimates on norms.\\
The following properties of $\Vert\cdot\Vert$ are immediate consequences of definition:
\beq\label{inv}1\leq \Vert x\Vert=\Vert x^{-1}\Vert, \quad x\in G,\eeq
\beq\label{prod} \Vert xy\Vert \leq \Vert x\Vert \Vert y\Vert,\quad x, y\in G.\eeq

In order to have estimates, we introduce the following notation.
Let $r$ be a positive integer. Let $f$ and $g$ be two positive functions defined over a subset $W$ of $G^r$ . 
\ber\label{comparfg} We write $f(x)\preccurlyeq  g(x), x\in W $ if and only if there are a positive constant $c$ and a positive integer $d$ such that $f(x)\leq c   g(x)^{d}$   for all $x\in W$.\eer
\ber\label{equival} We write  $f(x)\approx g(x), x\in\ W$ if   $f(x)\preccurlyeq g(x), x\in W  $ and $g(x)\preccurlyeq  f(x), x\in W $.\eer\\
If $f_1, f_2$ and $f_3$ are positive functions on $G^r$, we clearly have\me

if $f_1(x)\preccurlyeq f_2(x), x\in W$ and $f_2(x)\preccurlyeq f_3(x), x\in W$ then $f_1(x)\preccurlyeq f_3(x), x\in W$,

if $f_1(x)\approx f_2(x), x\in W$ and $f_2(x)\approx f_3(x), x\in W$ then $f_1(x)\approx f_3(x), x\in W$.\me

\no Moreover, if $f_1,f_2, g_1$ and $g_2$ are positive functions on $G^r$ which take values greater or equal to $1$, we obtain easily the following properties:
\ber\label{compar>1}\begin{enumerate}
\item for all positive integer $d$, we have $f_1(x)\approx f_1(x)^d, x\in W$,

\item if $f_1(x)\preccurlyeq g_1(x), x\in W$ and $f_2(x)\preccurlyeq g_2(x), x\in W$ then\\ $(f_1f_2)(x)\preccurlyeq (g_1g_2)(x), x\in W$,
\item if $f_1(x)\approx g_1(x), x\in W$ and $f_2(x)\approx g_2(x), x\in W$ then\\ $(f_1f_2)(x)\approx (g_1g_2)(x), x\in W$.
\end{enumerate}\eer
\no Since $\Vert x\Vert=\Vert x y y^{-1}\Vert\leq \Vert xy\Vert \Vert y\Vert$ and $\Vert xy\Vert\leq \Vert x\Vert  \Vert y\Vert$, we obtain
 \ber\label{compact} If $\Om$ is  a compact subset of $G$, then $\Vert x\Vert\approx \Vert x\om\Vert, \quad x\in G, \om\in \Om.$\eer
 Let  $P=M_PN_P$ be a parabolic subgroup of $G$ containing $M_0$. Then, each $x\in G$ can be written $x=m_P(x)n_P(x)k$ where $m_P(x)\in M_P, n_P(x)\in N_P$ and $k\in K$.  By (\cite{W2} Lemma II.3.1), we have
\beq\label{mn} \Vert m_P(x)\Vert+\Vert n_P(x)\Vert \preccurlyeq \Vert x\Vert, \quad x\in G.\eeq\\
Recall that $G^1$ is the   kernel of $h_{G}: G\to a_{G}$. Let us prove that 
\beq\label{g1a} \Vert xa\Vert\approx \Vert x\Vert \Vert a\Vert,\quad x\in G^1, a\in A_G.\eeq\\
According to  the Cartan decomposition (\ref{Cartan}), if $g\in G$, we denote by $m_0(g)$
an element of $M_0$ such that there exist $k,k'\in K$ with $g=km_0(g)k'$.  Notice that $\Vert h_{M_0}(m_0(g))\Vert$ does not depend on our choice of 
$m_0(g)$.  By (\ref{compact}), one has  
\beq\label{Normegm}\Vert g\Vert\approx\Vert m_0(g)\Vert,\quad g\in G,\eeq
 and by (\cite{W2}) 1.1.(6)) we have
\beq\label{Normem0} \Vert m_0 \Vert \approx e^{\Vert h_{M_0}(m_0 )\Vert}, \quad m_0\in M_0.\eeq
Let $x\in G^1$ and $a\in A_G$. Then $m_0(x)\in G^1\cap M_0$ and  $m_0(xa)=m_0(x) a$. Thus, one has  $h_G(m_0(x))=0$. We deduce from  (\ref{hGM}) 
 that 
$h_{M_0}(m_0(x))$ belongs to $a_{0}^G.$
Since $ h_{M_0}(m_0(x)a)=h_{M_0}(m_0(x))+h_{M_0}(a)$ and $h_{M_0}(a)\in a_G$,
we obtain by orthogonality that
$$\frac{1}{2}(\Vert h_{M_0}(m_0(x))\Vert +\Vert h_{M_0}(a)\Vert)\leq \Vert h_{M_0}(m_0(x)a)\Vert \leq \Vert h_{M_0}(m_0(x))\Vert+ \Vert h_{M_0}(a)\Vert.$$
\no Hence  (\ref{g1a}) follows from (\ref{Normegm}) and (\ref{Normem0}).\me

We denote by $C_c^\infty(G)$ the space of smooth functions on $G$ with compact support.
We normalize Haar measures according to  \cite{ArLT} \textsection 1. Unless otherwise stated, the Haar measure on a compact group will be normalized to have total volume $1$.

Let $M$ be a Levi subgroup of $G$.  
We fix  a Haar measure on $a_M$ so that the volume of the quotient $a_M/\tilde{a}_{M,\F}$ equals $1$.   

Let   $P=M N_P\in\cP(M)$. We denote by $\de_P$ the modular function of $P$   given by
$$\de_P(mn)=e^{2\rho_P(h_M(m))}, m\in M, n\in N_P,$$
 where $2\rho_P$ is the sum of roots, with multiplicity, of $(P,A_M)$. Let  $\bar{P}=M N_{\bar{P}}$ be the  the parabolic subgroup which is opposite to $P$. If $dn$   is a Haar measure on $N_P$ then the number $$\ga(P)=\int_{N_P} e^{2\rho_{\bar{P}}(h_{\bar{P}}(n))} dn$$ is finite. Moreover, the measure  $\ga	(P)^{-1} dn$ is independent of the choice of $dn$ and thus defines a canonical Haar measure on $N_P$. \\
 If $dm$ is a Haar measure on $M$ then there exists a unique Haar measure $dg$ on $G$, independent  of the choice of the parabolic subgroup   $P$,  such that
 $$\int_G f(g) dg= \frac{1}{\ga(P)\ga(\bar{P})}\int_{N_P}\int_M\int_{N_{\bar{P}}} f(nm\bar{n}) \de_P(m)^{-1} d\bar{n}\; dm\; dn,$$
 for 
$f\in C_c^\infty(G)$.   We say that $dm$ and $dg$ are compatible. Compatibility has the obvious transitivity property relative to Levi subgroups of $M$. Using the Iwasawa decomposition (\ref{Iwasawa}), these measures satisfy
$$\int_G f(g) dg=\frac{1}{\ga(P)}\int_K\int_M\int_{N_P} f(mnk) dn\; dm \; dk.$$

\subsection{The symmetric space $H\bb G$}
 Let   $\E$ be an unramified  quadratic extension of  $ \F$. Thus  $ \E=\F[\tau]$ where $\tau^2 $  is not a square in  $ \F$. We denote by $\si$ the  nontrivial element of the Galois group $\cG al(\E/\F)$ of $\E/\F$. The normalized valuation $|\cdot|_\E$ on $\E$ satisfies $|x|_\E=|x|_\F^2$ for $x\in\F$.\me

If $\sJ$ is an algebraic group defined over $\F$, as usual we denote  by $J$ its group of points over $\F$. Let  $\sJ\times_\F\E$ be the group, defined over $\E$,  obtained from $\sJ$ by extension of scalars. We consider the group 
$$\underline{\tilde{J}}:=\mbox{Res}_{\E/\F}(\sJ\times_\F\E)$$
defined over $\F$, obtained by   restriction of scalars.\me

With our convention, one has $\tilde{J}=\underline{\tilde{J}}(\F)$ and $\tilde{J}$ is isomorphic to $\sJ(\E)$.\me

 Let $\sH$ be a reductive group defined over $\F$. In all this article, we assume that $\sH$ is split over $\F$ and we set $\sG:=\underline{\tilde{H}}$ and $G:=\tilde{H}$. We fix a maximal split torus $A_0$ of $H$. Then $A_0$ is also a maximal split torus of $G$ and we have $A_H=A_G$. \me
 

The nontrivial element $\si$ of $\cG al(\E/\F)$ 	induces  an involution of $\sG$ defined over $\F$, which we denote by the same letter. This automorphism $\si$ extends to an $\E$-automorphism $\si_\E$ on $\sG\times_\F\E$. 

We consider  the canonical map $\varphi$ defined over $\F$ from $\sG$ to $(\sH\times_\F\E)\times (\sH\times_\F\E)$  by $\varphi (g)=(g,\si(g))$. \ber\label{isoProduit}  Then,   $\varphi$ extends uniquely to  an isomorphism $\Psi$ defined over $\E$ from $\sG\times_\F\E$ to $(\sH\times_\F\E)\times (\sH\times_\F\E)$ such that $\Psi(g)=(g,\si(g))$  for all $g\in \sG $ and if $\Psi(g)=(g_1,g_2)$ then $\Psi(\si_\E(g))=(g_2,g_1)$. 
\eer
Now we turn to the description of the geometric structure of the symmetric space $\cS=H\bb G$ according to \cite{RR}  sections 2 and 3.\me

Let  $\underline{\g g}$ be the  Lie algebra of $\sG$ and  $\g g$ be the Lie algebra of its  $\F$-points. We will say that  $\g g$  is the Lie algebra of  $G$ and the Lie algebra $\g h$ of $H$ consists of the elements of $\g g$ invariant by $\si$. We denote by  $\g q$ the space of antiinvariant elements of $\g g$ by $\si$. Thus, 
one has  $\g g=\g h\oplus \g q$ and  $\g g$ may be identified with $\g h\otimes_\F \E$.

As in  (\cite{RR} \textsection 2.), we say that a subspace $\g c$ of $\g q$ is a Cartan subspace of $\g q$ if $\g c$  is a maximal abelian subspace of $\g q$ made of semisimple elements. As $\E=\F[\tau]$, the multiplication by $\tau $ induces an isomorphism between the set of Cartan subspaces of $\g q$ and the set of Cartan subalgebras of $\g h$ which preserves $H$-conjugacy classes. \me

We denote by  $\scP$  the connected component of $1$ in the set of $x$ in $\sG$ such that $\si(x)=x^{-1}$. Then the map $\sp$ from $\sG$ to $\scP$ defined by $\sp(x)=x^{-1}\si(x)$ induces  an isomorphim of affine varieties $\sp: \sH\bb \sG\to \scP$.\me

A torus $\sA$ of $\sG$ is called a $\si$-torus if $\sA$ is a torus defined over $\F$ contained in $\scP$. Notice that such torus are called $\si$-split torus in \cite{RR}. We prefer change the terminology as $\si$-tori are not necessarily split over $\F$. Each $\si$-torus is the centralizer in $\scP$ of a Cartan subspace of $\g q$, or equivalently of a Cartan subalgebra of $\g h$.
\me

Let $S$ be a maximal torus of $H$.  We denote by     $\sS_\sigma$ the connected component of $\widetilde{\sS}\cap \scP$. Then $\sS_\si$ is a $\si$-torus defined over $\F$ which identifies with the antidiagonal $\{(s, s^{-1}); s\in \sS\}$ of $\sS\times\sS$ by the isomorphism (\ref{isoProduit}). Thus, $\sS_\si$ is a maximal $\si$-torus  and each maximal $\si$-torus arises in this way . The  $H$-conjugacy classes of   maximal tori of $H$ are in bijective correspondence with the $H$-conjugacy classes  of maximal $\si$-tori of $G$ by the map $S\mapsto S_\si$.
The roots  of $\sS$ (resp.; $\sS_\si$) in $\underline{\g h}=Lie(\sH)$ (resp.; $\underline{\g q}\otimes_\F\bar{\F}$) are the restrictions of the roots of $\tilde{\sS}$ in $\underline{\g g}=Lie(\sG)$. 
\ber\label{roots}Therefore, each root of $\sS$ (resp.; $\sS_\si$) in $\underline{\g g}$ has multiplicity two. If $\tilde{\sS}$ splits over a finite extension $\F'$  of $\F$, we denote by $\Phi ( S_\si', \g g')$ (resp.; $\Phi(S',\g h')$) the set of  roots of $\sS_\si(\F')$ in $ \g g \otimes_\F\F'$ (resp.; $\sS(\F')$ in $ \g h \otimes_\F\F'$).\\
Let $\tilde{\underline{\g s}}$ be the Lie algebra of $\tilde{\sS}$. Then, the differential of  each root $\al$ of $ \Phi (\tilde{S}',\g g')$ defines a linear form on $\tilde{\g s}\otimes_\F\F'$ which we denote by the same letter.
\eer

Let $\cG al(\overline{\F}/\F)$ be the Galois group of $\overline{\F}/\F$. By (\cite{RR} \textsection  3), the set of $(H,S_\si)$-double cosets in $\sH\sS_\si\cap G$ are parametrized by the finite set  $I$ of cohomology classes in $H^1(\cG al(\overline{\F}/\F), \sH\cap \sS_\si)$ which split in both $\sH$ and $\sS_\si$. To each such classe $m$,  we attach an element $x_m\in G$ of the form $x_m=h_m a_m^{-1}$ with $h_m\in\sH$ and $a_m\in\sS_\si$ such that $m_\ga=h_m^{-1} \ga(h_m)=a_m^{-1}\ga(a_m)$ for all $\ga\in\cG al(\overline{\F}/\F)$. 

\begin{lem}\label{fixeAS} Let $x\in  G$ such that $x=hs$ with $h\in\sH$ and $s\in \tilde{\sS}$. Then, $xSx^{-1}$ is a maximal torus of $H$ and there exists $h'\in H$ such that $x'=h'x$ centralizes the split connected component $A_S$  of $S$. 
\end{lem}
\dem Replace $S$ by a $H$-conjugate   if necessary, we may assume that $A:=A_S$ is contained in the fixed maximal split torus  $A_0$  of $H$.  Since $H$ is split,  $A_0$ is also a maximal split torus of $G$. 
	
Since  $x=hs\in G$, the torus $\sS':=x\sS x^{-1}$ is equal to $h\sS h^{-1}\subset \sH$. Thus $ \sS'$ is    defined over $\F$ and contained in $\sH$ and we obtain the first assertion.\me

Let $S':=\sS'(\F)$ and let $A'$ be the split connected component of $S'$. There exists $h_1\in H$ such that $h_1A'h_1^{-1}\subset A_0$. We set $x_1=h_1 x$, thus we have $A_1:=x_1 A x_1^{-1}\subset A_0$.
 
 Let $M=Z_G(A)$ and $M_1=Z_G(A_1)= x_1 M x_1^{-1}$. Then $A_0$ and $x_1 A_0 x_1^{-1}$ are maximal split tori of $M_1$. Therefore, there exists $y_1\in M_1$ such that $y_1 x_1 A_0 x_1^{-1} y_1^{-1}=A_0$. As $H$ is split, the Weyl group of $A_0$ in $G$ coincides with the Weyl group of $A_0$ in $H$. Thus, there exist $h_2\in N_H(A_0)$ and $v\in Z_G(A_0)$ such that   $z:=y_1 x_1= h_2 v $. 
 
 For $a\in A\subset A_0$, one has $z a z^{-1}= h_2 a h_2^{-1}= y_1 x_1a x_1^{-1} y_1^{-1}= x_1a x_1^{-1}$ since $ x_1a x_1^{-1}\in A_1$ and $y_1\in M_1$. One deduces that $x':=h_2^{-1}h_1x$ centralizes $A$.\qed
 \ber\label{xm}Thus, for  each maximal torus $S$ of $H$, we can fix a finite set of representatives  $\kappa_S=\{x_m\}_{m\in I}$ of the  $(H,S_\si)$-double cosets in $\sH\sS_\si\cap G$ such that each element $x_m$  may be written   $x_m=h_m a^{-1}_m$ where  $h_m\in\sH$  centralizes   $A_S$ and $a_m\in \sS_\si$. Hence $x_m$ centralizes $A_S$.\eer
\subsection{Weyl integration formula and orbital integrals}

We first recall   basic notions on the symmetric space  according to (\cite{RR}, \textsection 3). An element $x$ in $\sG$ is called $\si$-semisimple if the double coset $\sH x\sH$ is Zariski closed. This is equivalent to say that $\sp(x)$ is a semisimple point of $\sG$. We say that a semisimple element $x$ is $\si$-regular if this closed double coset $\sH x\sH$  is of maximal dimension. This is equivalent to say that the centralizer of $\sp(x)$ in $\g q$ (resp.; $\scP$) is a Cartan subspace of $\g q$ (resp.; a maximal $\si$-torus of $\sG$).\me

We denote by $G^{\si-reg}$ the set of $\si$-regular elements of $G$.\me

For $g\in G$, we denote by  $D_G(g)$ the coefficient of the least power of $t$ appearing nontrivially in $\det(t+1-\textrm{Ad} (g))$. We define the $H$-biinvariant function $\Delta_\si$ on $G$ by 
$\Delta_{\si}(x)=D_G(\sp(x))$. Then by (\cite{RR}, Lemma 3.2. and Lemma 3.3), the set of $g\in G$ such that  $\Delta_\si(g)\neq 0$ coincides with   $G^{\si-reg}$. \\
Let $S$ be a maximal torus of $H$ with Lie algebra $\g s$. Then $\tilde{\g s}:=\g s\otimes_\F\E$ identifies with the Lie algebra of $\tilde{S}$. For $g\in x_m S_\si$ with $x_m\in \kappa_S$, one has
\beq\label{deltas} \Delta_\si(g)=D_G(\sp(g))=\det(1-\textrm{Ad}(\sp(g)))_{\g g/\tilde{\g s}}.\eeq

\no By (\cite{RR}  Theorem 3.4 (1)), the set  $G^{\si-reg}$ is a disjoint union
\ber\label{Thm3.4RR}
 $$ G^{\si-reg}=\bigcup_{\{S\}_H}\bigcup_{x_m\in \kappa_S} H\big( (x_mS_\si)\cap G^{\si-reg}\big)H,$$\\
 where $\{S\}_H$ runs the $H$-conjugacy classes of   maximal tori of $H$.\eer\\
 If $x_m\in\kappa_S$ then $x_m=h_ma_m$ for some $h_m\in\sH$ and $a_m\in\sS_\si$, hence $\sp(x_m)=a_m^{-2}$ commutes with $S$ and $S_\si$. Therefore for $\ga\in S_\si$, we have
 $$\sp(x_m\ga)=\sp(x_m)\ga^{-2}\quad\mbox{and}\quad Hx_m\ga S=Hx_m\ga.$$\\
We have the following   Weyl integration formula (cf.  \cite{RR} Theorem 3.4 (2)):

\ber\label{Weyl0} Let $f$ be a compactly supported smooth function on $G$, then we have $$\int_G f(y) dy=\sum_{\{S\}_H}\sum_{x_m\in \kappa_S} c^0_{S,x_m} \int_{S_\si} |\De_\si(x_m \ga)|_\F^{1/2}\int_{S\bb H}\int_H f(hx_m\ga l)dh d\bar{l}d\ga,$$ \eer
where the constants $c^0_{S, x_m}$ are explicitly given in (\cite{RR} Theorem 3.4 (1)).\me

For our purposes, we need another version of this Weyl integration formula.\\
Let $S$ be a maximal torus of $H$. We denote by $A_S$  its split connected component.  Since the  quotient $A_S\bb S$ is  compact, by our choice of measure, the integration over $S\bb H$ in the Weyl formula above can be replaced by an integration over $A_S\bb H$. Moreover, it is convenient to change $h$ into $h^{-1}$. As every $x_m\in\kappa_S$ commutes with $A_S$ (cf.   (\ref{xm})), one can replace the integration over $(A_S\bb H) \times H$, by an integration over $diag(A_S)\bb (H\times H)$ where $diag(A_S)$ is the diagonal of $A_S$. This gives the following   Weyl integration formula equivalent to (\ref{Weyl0}):
\beq\label{Weyl1}\int_G f(y) dy=\sum_{\{S\}_H}\sum_{x_m\in \kappa_S} c^0_{S,x_m} \int_{S_\si} |\De_\si(x_m \ga)|_\F^{1/2}\int_{diag(A_S)\bb (H\times H)}f(h^{-1}x_m\ga l) d\overline{(h,l)}d\ga .\eeq

We will now describe the  $H$-conjugacy classes of   maximal tori of $H$ in terms of   Levi subgroups $M$ in $\cL(A_0)$ and $M$-conjugacy classes of some tori of $M$.

Let $M\in \cL(A_0)$. We denote by $N_H(M)$ its normalizer  in $H$.   If $S$ is a maximal torus of   $M$, we denote by $W(M,S)$ (resp. $W(H,S)$) its Weyl group in $M$ (resp. $H$). We choose a set $\cT_M$ of representatives for the $M$-conjugacy classes of   maximal tori $S$ in $M$ such that $A_M\bb S$ is compact. For $M,M'\in\cL(A_0)$, we write $M\sim M'$ if $M$ and $M'$ are conjugate under $H$. \me

Let $S$ be a maximal torus of $H$ whose split connected component $A_S$ is contained in $A_0$. Then, the centralizer $M$ of $A_S$ belongs to $\cL(A_0)$ and $S$ is a maximal torus of $M$ such that $A_M\bb S$ is compact. If  $S'$ is a maximal torus conjugated to $S$ by $H$ such that   $A_{S'} $ is contained in $A_0$, then the centralizer $M'$  of $A_{S'} $   in $H$ belongs to $\cL(A_0)$ and $M'\sim M$. 

Since each maximal torus of $H$ is $H$-conjugate to a maximal torus $S$ such that $A_S\subset A_0$, we obtain a surjective map $S\mapsto \{S\}_H$ from the set of  $S$ in $\cT_M$ where $ M$ runs a system of representatives of $\cL(A_0)_{/\sim}$ to the set of $H$-conjugacy classes of maximal tori of $H$. 

Let $M\in\cL(A_0)$.  By (\cite{Ko} (7.12.3)),  the cardinal of the class of $M$ in $\cL(A_0)_{/\sim}$ is equal to  $$\frac{\vert W(H,A_0)\vert }{\vert W(M,A_0)\vert \vert N_H(M)/M\vert}$$ where $N_H(M)$ is the normalizer of $M$ in $H$.\\
  By (\cite{Ko} Lemma 7.1), if $S$ is a maximal torus of $M$, then the  number of $M$-conjugacy classes of maximal torus $S'$ in $M$ such that $S'$ is $H$-conjugate to $S$ is equal to
$$\frac{\vert N_H(M)/M\vert \vert W(M,S)\vert }{\vert W(H,S)\vert}.$$ \\
Therefore, we can rewrite (\ref{Weyl1}) as follows:
 \beq\label{Weyl2} \int_G f(g) dg=\sum_{M\in\cL(A_0)} c_M\sum_{S\in\cT_M} \sum_{x_m\in \kappa_S}c_{S, x_m}\int_{S_\si} |\De_\si(x_m\ga)|_\F^{1/2} \int_{diag(A_M)\bb H\times H}  f(h^{-1} x_m\ga l)d\overline{(h,l)} d\ga \eeq
where $$c_M=\dfrac{|W(M, A_0)|}{|W(H,A_0)| }\quad \mbox{ and } \quad c_{S,x_m}=\dfrac{\vert W(H,S)\vert}{\vert W(M,S)\vert}c^0_{S,x_m}.$$\me

 Let $f\in C_c^\infty(G)$. We define the orbital integral $\cM(f)$  of $f$ on $G^{\si-reg}$ as follows. Let $S$ a maximal torus of $H$. For $x_m\in\kappa_S$ and $\ga\in S_\si$ with $x_m\ga\in G^{\si-reg}$, we set
\ber\label{OI}$$\cM(f)(x_m\ga):=|\De_\si(x_m\ga)|_\F^{1/4}\int_{ diag(A_S)\bb (H\times H)} f(h^{-1} x_m\ga l) d\overline{(h,l)}$$
$$=|\De_\si(x_m\ga)|_\F^{1/4}\int_{S\bb H}\int_H f(h x_m\ga l) dh d\overline{l}.$$\eer
Our definition corresponds, up to a positive constant, to Definition 3.8 of \cite{RR}.  Indeed, by definition of $\De_\si$,  we have 
$\De_\si(x_m\ga)=D_G(\sp(x_m\ga))$.
Since we can write $x_m=h_m a_m$ with $h_m\in\sH$ and $a_m\in\sS_\si$, we have $\sp(x_m\ga)=\sp(x_m)\ga^{-2}=a_m^{-2} \ga^{-2}$ for $\ga\in S_\si$. Let $\F'$ be an extension of $\E$ such that  $\tilde{S}$ splits over $\F'$ and $a_m\in \sS_\si(\F')$. Since each  root $\al$ of $\sS_\si(\F')$ in $\g g\otimes\F'$ have multiplicity $m(\al)=2$, using notation of (\ref{roots}), we obtain:
$$\De_\si(x_m\ga)=\prod_{\al\in \Phi(S'_\si,\g g')}(1-\sp(x_m)^\al\ga^{-2\al})^2=\prod_{\al\in \Phi(S'_\si,\g g')}(\ga^\al-\sp(x_m)^\al\ga^{-\al})^2,$$
hence $$\vert \De_\si(x_m\ga)\vert_{\F'}^{1/4}=\prod_{\al\in \Phi(S'_\si,\g g')}\vert(\ga^\al-\sp(x_m)^\al\ga^{-\al})^{m(\al)-1}\vert_{\F'}^{1/2},$$
$$=\prod_{\al\in \Phi(S'_\si,\g g')}\vert(\ga^\al-\sp(x_m)^\al\ga^{-\al})\vert_{\F'}^{1/2}.$$\\
Then, the Weyl integration formula (\ref{Weyl0})  in terms of orbital integrals is given as in (\cite{RR} page 126) by
 $$\int_G f(y) dy=\sum_{\{S\}_H}\sum_{x_m\in \kappa_S} c^0_{S,x_m} \int_{S_\si} |\De_\si(x_m \ga)|_\F^{1/4} \cM(f)(x_m\ga)d\ga.$$

 \begin{theo}\label{OIbounded}  Let $f\in C_c^\infty(G)$ and $S$ be a maximal torus of $H$. Let $x_m\in\kappa_S$.
 \begin{enumerate} \item There exists a compact set $\Omega$ in $S_\si$ such that, for any  $\ga$ in the complementary of $\Om$ with  $x_m\ga\in G^{\si-reg}$, one has $\cM(f)(x_m\ga)=0$.
 \item  $$\sup_{\ga\in S_\si;\; x_m\ga\in G^{\si-reg}} |\cM(f)(x_m\ga)|<+\infty.$$
 \end{enumerate}
 \end{theo}
 \dem
 The proof follows that of  the group case (cf.  \cite{HCvD}  proof of Theorem 14). We write it  for convenience of the reader.\\
{\it 1.}  Let $\om$ be the support of $f$. We consider the set $\om_S$ of elements $\ga$ in $S_\si$ such that $x_m\ga$ is in the closure of $H\om H$.  For $g\in G$, we consider the polynomial function
\beq\label{DG} \det (1-t-\Ad\; \sp(g))=t^n+q_{n-1}(g) t^{n-1}+\ldots + q_l(g) t^l\eeq
 where $l$ is the rank of $G$ and $n$ its dimension. Each $q_j$ is a $H\times H$ biinvariant  regular function on $G$, thus it is bounded on $x_m\om_S$. Therefore, the roots of $\det (1-t-\Ad\; \sp(g))$ are bounded on  $x_m\om_S$. 
 
For $\ga\in S_\si$, we have $\sp(x_m\ga)=\sp(x_m)\ga^{-2}$. We choose a  finite extension $\F'$  of $\F$ such that $\tilde{\sS}$ splits over $\F'$ and $\sp(x_m)\in\sS_\si(\F')$. Using notation of  (\ref{roots}),  the roots of $ \det (1-t-\Ad\; \sp(x_m\ga))$ are the numbers  $(1-\sp(x_m)^\al \ga^{-2\al})$ for $\al\in\Phi(S_\si', \g g')$. Since these roots are bounded on $x_m\om_S$, 
we deduce  that    the maps $\ga\to \ga^{\al}$, $\al\in \Phi(S_\si', \g g')$, are bounded on $\om_S$. This implies that  $\om_S$ is bounded. Then, the closure $\Om$ of $\om_S$ satisfies the first assertion.\me
 
 
\noindent {\it 2.}  By {\it 1.}, if $\ga\notin \Om$ then $\cM(f)(x_m\ga)=0$. Thus, it is enough to prove that for each $\ga_0\in S_\si$, there exists a neighborhood $V_{\ga_0}$ of $\ga_0$ in $S_\si$ such that 
\beq\label{majOI}\sup_{\ga\in V_{\ga_0},  x_m\ga\in G^{\si-reg}} |\cM(f)(x_m\ga)|<+\infty.\eeq

\no  Let $y_0:=\sp(x_m\ga_0)$. We first  
assume that $y_0$ is central in $G$. Then, we have  $\De_\si(x_m\ga_0\ga)=D_G(y_0\ga^{-2})=D_G(\ga^{-2})$ for $\ga\in S_\si$ and  $x_m\ga_0h (x_m\ga_0)^{-1}\in H$ for $h\in H$. We define the function $f_0$ on $G$ by $f_0(g):=f(x_m\ga_0g)$. Then, we have $\cM(f_0)(\ga)=\cM(f)(x_m\ga_0\ga)$ for $\ga\in S_\si\cap G^{\si-reg}$. Thus, we are reduced  to the case $y_0=1$. As in the group  case, we use the exponential map $"\exp"$ which is well defined in a neighborhood of $0$ in $\g g$ since the characteristic of $\F$ is equal to zero (cf.  \cite{HCBS}  \textsection10). As in (\cite{HC}  proof of Lemma 15), we can choose a $H$-invariant open   neighborhood $V_0$ of $0$ in $\g h$ such that the map $X\in V_0\mapsto \exp (\tau X)$ is an isomorphism and an homeomorphism onto its image and there is a $H$-invariant function $\varphi\in C^\infty_c(\g h)$ such that $\varphi(X)=1$ for $X\in V_0$. We define $\bar{f}$ in $C_c^\infty(\g h)$ by $\bar{f}(X)= \varphi (X) \int_H f(h \exp (\tau X)) dh$.
 
 Let $\g s$ be the Lie algebra of $S$. For $X\in \g s$, we set $\eta(X)=|\det (\textrm{ad} X)_{\g h/\g s}|_\F$. We consider a finite extension $\F'$ of $\F$  such that $\tilde{\sS}$  splits over $\F'$ and $\sp(x_m)\in \sS_\si(\F')$. We use   notation  of (\ref{roots}). Since each root of $S'_\si$ in $\g g'$ has  multiplicity $2$, we have for $X\in V_0$
 $$\frac{|\De_\si(\exp\tau X)|_{\F'}^{1/2}}{\eta(X)}=\frac{|D_{G'}(\exp(-2\tau X)|_{\F'}^{1/2}}{\eta(X)}=  \frac{\prod_{\al\in\Phi( S' ,\g h')}\vert 1-e^{2\tau\al(X)}\vert_{\F'}}{\prod_{\al\in\Phi( S' ,\g h')}\vert \al(X)\vert_{\F'}}$$
 $$=|2\tau|_{\F'}^{|\Phi( S' ,\g h')|}\prod_{\al\in \Phi( S' ,\g h')}|1+ \tau\al(X) +\frac{4\tau^2\al(X)^2}{3!}+\ldots |_{\F'}.$$
 We can reduce $V_0$ in such way that each term of this product is equal to $1$. Thus, we obtain
 $$\cM(f)(\exp\tau X)= |2\tau|_{\F'}^{|\Phi( S' ,\g h')|} \eta(X)^{1/2}\int_{H/S}\big(\int_H f(h\exp \tau \Ad(l)X) dh\big)d\bar{l}$$
 $$=  |2\tau|_{\F'}^{|\Phi( S' ,\g h')|} \eta(X)^{1/2}\int_{H/S}\bar{f}(\Ad(l)X) d\bar{l},$$
  for $X\in V_0$.  The estimate (\ref{majOI}) follows from the result on the Lie algebra (cf. \cite{HCvD}  Theorem 13).\me

 If $y_0=\sp(x_m\ga_0)$ is not central in $G$, we consider the centralizer $\scZ$ of $y_0$ in $\sH$. Let $\scZ^0$ be the connected neutral component of $\scZ$. By(\cite{Bo}, III.9), the group $\scZ^0$ is defined over $\F$. As usual, we set $\tilde{\underline{\cZ}}^0:=\textrm{Res}_{\E/\F}(\scZ^0\times_\F\E)$ and we denote by  $\tilde{\underline{\g z}}$ its Lie algebra.  By definition of $\tilde{\g z}$, one has $$|\det (1-\textrm{Ad}(y_0))_{\g g/\tilde{\g z}}|_\F\neq 0.$$ 
 \no Thus, there exists  a neighborhood $V$ of $1$ in $S_\si$ such that, for all $\ga\in V$, then 
 \beq\label{det0}|\det (1-\textrm{Ad}(y_0\ga^{-2}))_{\g g/\tilde{\g z}}|_\F=|\det (1-\textrm{Ad}(y_0)_{\g g/\tilde{\g z}})|_\F\neq 0.\eeq

From (\cite{HCvD}  Lemma 19), there exist a neighborhood $V_1$ of $y_0$ in $\tilde{S}$ and  a compact subset $\overline{C_G}$ of $\tilde{\cZ}^0\bb G$ such that, if $g\in G$ satsifies $g^{-1} V_1 g\cap \sp(\om)\neq \emptyset$ then its image $\bar{g}$ in $\tilde{Z}^0\bb G$ belongs to  $ \overline{C_G}$  (here $\om$ is the support of $f$). 

We choose a neighborhood $W$ of $1$ in $S_\si$ such that $W\subset V$ and $\sp(x_m\ga_0\ga)=y_0 \ga^{-2}\in V_1$ for all $\ga\in W$. By (\cite{Bo}, III 9.1), the quotient $\cZ^0\bb H$ is a closed subset of  $\tilde{\cZ}^0\bb G$, hence 
\ber\label{Cbar}the set $\overline{C}:=\overline{C_G}\cap \cZ^0\bb H$ is a compact subset of $\cZ^0\bb H$ such that
if $l\in H$ satisfies $l^{-1}y_0\ga^{-2}l\in\sp(\om)$ for some $\ga\in W$ then its image $\bar{l}$ in $\cZ^0\bb H$ belongs to $ \overline{C}$.\eer\\
Let  $\ga\in W$ such that $x_m\ga_0\ga\in G^{\si-reg}$. One has
\beq\label{express1}\int_{S\bb H} \int_H f(hx_m\ga_0\ga l) dh d\bar{l}=\int_{\cZ^0\bb H}\int_{S\bb \cZ^0} \int_Hf(hx_m\ga_0\ga\xi l)dh  d\bar{\xi} d\bar{l}.\eeq
By the choice of $W$, the map  $$ \bar{l}\in \cZ^0\bb H \mapsto \int_{S\bb \cZ^0} \int_Hf(hx_m\ga_0\ga\xi  l)dh  d\bar{\xi} $$
vanishes outside $\bar{C}$.
We choose $u\in C_c^\infty(H)$ such that   the map $\overline{u}\in C_c^\infty(\cZ^0\bb H)$ defined by  $\overline{u}(\bar{l}):=\int_{\cZ^0}u(\xi l) d\xi$ is equal to $1$ if $\bar{l}\in\overline{C}$. As $u$ and $f$ are  compactly supported, the map
$$\Phi: z\in\tilde{\cZ}^0\mapsto  \int_H u(l) \int_H f(h x_m\ga_0z l) dh dl$$
is well-defined. Since $y_0=\sp(x_m\ga_0)=(x_m\ga_0)^{-1}\si(x_m\ga_0)$, we have $\xi (x_m\ga_0)^{-1}\si(x_m\ga_0)=(x_m\ga_0)^{-1}\si(x_m\ga_0)\xi$ for $\xi\in \cZ^0$. Hence, $x_m\ga_0 \xi(x_m\ga_0)^{-1}\in H$. Thus $\Phi$  is left invariant by $\cZ^0$. 

We claim that $\Phi\in C_c^\infty(\cZ^0\bb \tilde{\cZ^0})$. Indeed, fix $l$ in the support of $u$. If $f(h x_m\ga_0z l)$ is nonzero for some  $h\in H$ and $z\in\tilde{\cZ}^0$ then $\sp(h x_m\ga_0z l)=\sp( x_m\ga_0z l)$ belongs to $\sp(\om)$, where $\om$ is the support of $f$. Since $z$ commutes with $y_0=\sp(x_m\ga_0)$, we have $\sp( x_m\ga_0z l)=l^{-1} y_0\sp(z)\si(l)$. As $u$ is compactly supported, we deduce that $\Phi(z)=0 $ when $\sp(z)$ is  outside a compact set. Hence, the map $\Phi$ is a compactly supported function on $\cZ^0\bb\tilde{\cZ}^0$.  

By assumption, the function  $f$ is right  invariant by a compact open subgroup  of $G$. Thus $f$ is right invariant by some compact open subgroup of $H$.   We denote by $\tau_lf$ the right translate of $f$ by an element $l\in G$. Since $u$ is compactly supported, the vector space generating by $\tau_lf$, when $l\in H$ runs the support of $u$, is finite dimensional. Hence, one can find a compact open subgroup $J_1$ of $\tilde{Z}^0$ such that for each $l$ in the support of $u$,  the function $\tau_lf$ is right invariant by $J_1$. This implies that $\Phi$ is smooth and our claim follows. 

 Therefore, there exists $\varphi\in C_c^\infty (\tilde{\cZ}^0)$ such that 
$$\Phi(z)=\int_{\cZ^0}\varphi(\xi z) d\xi= \int_H u(l) \int_H f(h x_m\ga_0z l) dh dl, \quad z\in\tilde{\cZ}^0.$$
We obtain 
$$\int_{S\bb \cZ^0} \int_{\cZ^0}\varphi(\xi_1 \ga \xi_2) d\xi_1 d\bar{\xi_2}= \int_H u(l)\big(\int_{S\bb \cZ^0}  \int_H f(h x_m\ga_0  \ga \xi_2 l) dh d\bar{\xi_2}\big) dl $$
$$= \int_{\cZ^0\bb H}\int_{\cZ^0} u(\xi_1l) \big(\int_{S\bb \cZ^0}\int_H f(h x_m\ga_0  \ga \xi_2 \xi_1l) dh d\bar{\xi_2}\big) d\xi_1d\bar{l}$$
 $$=\int_{\cZ^0\bb H} \overline{u}(\bar{l})  \big(\int_{S\bb \cZ^0}\int_H f(h x_m\ga_0  \ga \xi_2 l) dh d\bar{\xi_2}\big)  d\bar{l}.$$
 By  definition, the map $\bar{u}$ is equal to $1$ on the compact set $\overline{C}$. By definition of $\overline{C}$  (cf.  (\ref{Cbar}) and (\ref{express1})), we obtain
 $$\int_{S\bb \cZ^0} \int_{\cZ^0}\varphi(\xi_1 \ga \xi_2) d\xi_1 d\bar{\xi_2}= \int_{S\bb H} \int_H f(hx_m\ga_0\ga l) dh d\bar{l}.$$
By (\ref{det0}) and the choice of $W$,  one has $$|D_G(y_0\ga^{-2})|_\F=|D_{\tilde{\cZ}^0}(\ga^{-2})|_\F |\det (1-\textrm{Ad}(y_0))_{\g g/\tilde{\g z}}|_\F,\quad  \ga\in W .$$ Then, we deduce that for $\ga\in W$ satisfying $x_m\ga_0\ga\in G^{\si-reg}$, one has   $$\cM(f)(x_m\ga_0\ga)= |\det (1-\textrm{Ad}(y_0))_{\g g/\tilde{\g z}}|_\F^{1/4} |D_{\tilde{\cZ}^0}(\ga^{-2})|_\F^{1/4} \int_{S\bb \cZ^0} \int_{\cZ^0}\varphi(\xi_1 \ga \xi_2) d\xi_1 d\bar{\xi_2}.$$
Since $|D_{\tilde{\cZ}^0}(\ga^{-2})|_\F$ coincides with the function $\vert\De_\si\vert_\F$ for the group $\tilde{Z}^0$ evaluated at $\ga$ (cf. (\ref{deltas})),  one deduces the estimate (\ref{majOI}) for $f$ applying the first case to $\varphi$ defined on $\tilde{\cZ}^0$. 
\qed

\section{Geometric side of the local relative trace formula}

 \subsection{ Truncation}

In this section, we will recall some results of (\cite{ArLT}, \textsection 3), needed in the sequel. We keep notation of \textsection \ref{defgroups} for the group $H$. Since $H$ is split, one has $M_0=A_0$.
We fix a Levi subgroup $M\in\cL(A_0)$   of $H$. Let $P\in \cP(M)$. We recall that  $A_M$ denotes the maximal split connected component of $M$. 

We denote by $\Si_P$ the  set of roots of $A_M$ in the Lie algebra of $P$, $\Si_P^r$ the subset of reduced roots and $\De_P$ the subset of simple roots.  

For  $\be\in\De_P$,  the "co-root" $\check{\be}\in a_M$  is defined as usual  as follows: if $P\in\cP(A_0)$ is  a minimal parabolic subgroup, then $\check{\be}=2(\be,\be)^{-1}\be$, where $a_0^*$ identifies with $a_0$ by the scalar product on $a_0$. In the general case, we choose $P_0\in\cP(A_0)$ contained in $P$. Then, there exists a unique $\al\in\De_{P_0}$ such that $\be=\al_{|a_M}$. The "co-root" $\check{\be}$ is the projection of $\check{\al}$ onto $a_M$ with respect to the decomposition $a_0=a_M\oplus a_0^M$. This projection does not depend of the choice of $P_0$.\me

We denote by $a_P^+$ the positive Weyl chamber of elements $X\in a_M$ satisfying $\al(X)>0$ for all $\al\in \Si_P$.\me

Let $M\in \cL(A_0)$. A set of points in $a_M$ indexed by $P\in \cP(M)$
$$\cY=\cY_M:=\{Y_P\in a_M; P\in \cP(M)\}$$
is called a $(H,M)$-orthogonal set if  for all adjacent parabolic subgroups $P,P'$ in $\cP(M)$ whose chambers in $a_M$  share the wall determined by the simple root $\al\in\De_P\cap(-\De_{P'})$, one has $Y_P-Y_P'=r_{P,P'} \check{\al}$ for a real number $r_{P,P'}$. The orthogonal set is called positive if each of the numbers $r_{P,P'}$ are nonnegative. This is the case  for example if the number 
\beq\label{dY} d(\cY)=\inf_{\{\al\in\De_P; P\in\cP(M)\}} \al(Y_P)\eeq
is nonnegative.\\
One example is the set 
$$\{ -h_P(x); P\in\cP(M)\},$$
defined for any point $x\in H$. This is a positive $(H,M)$-orthogonal set by (\cite{ArDS} Lemma 3.6).
\ber\label{restY} If $L$ belongs to $\cL(M)$ and $Q$ is a group in $\cP(L)$, we define $Y_Q$ to be the projection onto $a_L$ of any point $Y_P$, with $P\in \cP(M)$ and $P\subset Q$. Then $Y_Q$ is independent of $P$ and $\cY_L:=\{Y_Q; Q\in \cP(L)\}$ is a $(H,L)$-orthogonal set.\eer
We shall write $\cS_M(\cY)$ for the convex hull in $a_M/a_H$ of a $(H,M)$-orthogonal set $\cY$. Notice that $\cS_M(\cY)$ does only depend on the projection onto $a_M^H$ of each $Y_P\in \cY$, $P\in\cP(M)$.\me

Let $P\in\cP(M)$. If each $Y_P$ is in the positive Weyl chamber $a_P^+$ (this condition is equivalent to say that   $d(\cY)$ is positive), we have a simple description of $\cS_M(\cY)\cap a_P^+$ (\cite{ArLT} Lemma 3.1). We denote by $(\om_\ga^P)_{\ga\in\De_P}$ the set of weights, that is the  dual basis in   $(a_M^H)^*$ of the set of co-roots $\{\check{\ga};\ga\in\De_P\}$. Then, we have

\beq\label{SMY-a+} \cS_M(\cY)=\{ X\in a_P^+; \om_\ga^P(X-Y_P)\leq 0, \ga\in \De_P\}.
\eeq
We now     recall a decomposition of the characteristic function of $\cS_M(\cY)$ valid when $\cY$ is positive. (cf.  \cite{ArLT} (3.8)). Suppose that $\La$ is a point in $a_{M,\C}^*$ whose real part $\La_R\in a_M^*$ is in general position. If $P\in \cP(M)$, we define $\De_P^\La$ the set of simple roots $\al\in \De_P$ such that $\La_R(\check{\al})<0$. Let $\varphi_P^\La$ be the characteristic function of the set of $X\in a_M$ such that $\om^P_\al (X)>0$ for each $\al\in \De_P^\La$ and $\om^P_\al(X)\leq 0$  for each $\al$ in the complement of $\De_P^\La$ in $\De_P$. We define
\beq\label{SigmaM} \si_M(X, \cY):= \sum_{P\in\cP(M)} (-1)^{|\De_P^\La|}\varphi_P^\La(X-Y_P).\eeq
\ber\label{SigmaMnul}By (\cite{ArLT}, \textsection 3 p22), the function $\si_M(\cdot, \cY)$  vanishes on the complement of $\cS_M(\cY)$ and is bounded. Moreover, if $\cY$ is positive then $\si_M(\cdot,\cY)$ is exactly the characteristic function of $\cS_M(\cY)$.\eer

The following Lemma will allow us to   define the minimum of two orthogonal sets.

For $P\in \cP(M)$, we denote by $(\tilde{\om}^P_\ga)_{\ga\in\De_P}$ the set of coweights, that is the dual basis in $a_M^H$ of the roots $\{\ga; \ga\in\in \De_P\}$.

\begin{lem}\label{copoids} Let $P$ and $P'$ two adjacent parabolic subgroups in $\cP(M)$ whose chambers in $a_M$  share the wall determined by the simple root $\al\in\De_P\cap(-\De_{P'})$. Then:
\begin{enumerate}\item For all $\be$ in $\De_P-\{\al\}$, there exists a unique $\be'$ in $\De_{P'}-\{\al\}$ such that $\be'=\be+k_\be\al$ where $k_\be$ is a nonnegative integer. Moreover, the map $\be\mapsto \be'$ is a bijection between $\De_P-\{\al\}$ and $\De_{P'}-\{-\al\}$.
\item  For all $\be$ in $\De_P-\{\al\}$, one has $\tilde{\om}^{P'}_{\be'}=\tilde{\om}_\be^P$.
\end{enumerate}\end{lem}
\dem We denote by $\N$ the set of nonnegative integers and by $\N^*$ the subset of positive integers.\\
{\it 1.} As $P$ and $P'$ are adjacent, we   have $\Si_{P'}=\big(\Si_P-\{\al\}\big)\cup\{-\al\}$.
Let $\be\in \De_P-\{\al\}$.  If $\be\in\De_{P'}$ then we set $\be':=\be$. \\
Assume  that $\be$ is not in  $\De_{P'}$. Since $\be\in\Si_{P'}$,  there exists $\Te\subset \De_{P'}-\{-\al\}$ such that   $\be=\sum_{\de\in\Te} n_\de\de-k_\be\al$  where the $n_\de$'s are positive integers and $k_\be$ is a nonnegative integer. Each $\de$ in $\Te$ belongs to $\Si_P$. Therefore, there are nonnegative integers $(r_{\de,\eta})_{\eta\in\De_P}$ such that 
$\de=\sum_{\eta\in\De_{P}}r_{\de,\eta} \eta$.  We set $\be_1:=\sum_{\de\in\Te} n_\de\de=\be+k_\be\al$.

Let $\ga\in\De_P-\{\al\}$. If $\ga\neq\be$, one has $\be_1(\tilde{\om}_\ga^P)= \be(\tilde{\om}_\ga^P)=0$. Thus, for each   $\de\in\Te$, we have $r_{\de,\ga}=0$, hence  $\de=r_{\de,\be}\be+r_{\de,\al}\al$.

On the other hand, one has    $\be_1(\tilde{\om}_\be^P)=\be(\tilde{\om}_\be^P)=1$. Thus, for all $\de\in\Te$, one has $\sum_{\de\in\Te} n_\de r_{\de,\be}=1$. Since $n_\de\in \N^*$ and $r_{\de,\be}\in \N$, one deduces that there exists a unique $\de_0\in\Te$ such that $r_{\de_0,\be}\neq 0$ and we have $n_{\de_0}=r_{\de_0,\be}=1$. This implies that $\Te=\{\de_0\}$ and $\be=\de_0-k_\be \al$.  We can take  $\be':=\de_0$. Hence, we obtain the existence of $\be'$ in all cases.\me

If $\be'_1\in \De_{P'}$ satisfies $\be'_1=\be+k^1_\be \al$ then $\be'= \be_1'+(k_\be-k^1_\be)\al$. Since the roots $\be'_1$, $\be'$ and $-\al$ belong to the set of simple roots $\De_{P'}$, we deduce that $\be_1'=\be'$. This gives the unicity of $\be'$. 

Let $\ga$ and $\be$ be in $\De_P$ such that $\ga'=\be'$. Then we have $\be=\ga+(k_{\ga}-k_\be)\al$.  Since $\ga,\be$ and $\al$ belong to $\De_P$,  the same argument as above leads to $\be=\ga$. Hence, the map $\be\mapsto \be'$ is injective.\me

\no {\it 2.} Let $\be\in \De_P-\{\al\}$. By  definition, we have $\be'=\be+k_\be\al\in \De_{P'}-\{-\al\}$ with $k_\be\in\N$. Thus  we have $\al(\tilde{\om}_{\be'}^{P'})=\al(\tilde{\om}_\be^P)=0$ and  $\be(\tilde{\om}_{\be'}^{P'})=\be'(\tilde{\om}_{\be'}^{P'})=1$.
If $\ga\in\De_P-\{\be,\al\}$, then $\ga'=\ga+k_\ga\al$ is different from $\be'$ by ({\it 1.}),  thus  we have $\ga(\tilde{\om}_{\be'}^{P'})=\ga'(\tilde{\om}_{\be'}^{P'})=0$. 
One deduces that $\tilde{\om}_{\be'}^{P'}=\tilde{\om}_{\be}^{P}$. \qed
\ber\label{infP} For $Y^1$ and $Y^2$ in $a_M$, we denote by $\inf^P\{Y^1,Y^2\}$ the unique element $Z$ in $a_M^H$ such that, for all $\ga\in \De_P$, one has $(\tilde{\om}_\ga^P, Z)=\inf\{ (\tilde{\om}_\ga^P, Y^1), (\tilde{\om}_\ga^P, Y^2)\}$.\eer
\begin{lem}\label{infY} Let $\cY^1=\{Y^1_P, P\in\cP(M)\}$ and  $\cY^2=\{Y^2_P, P\in\cP(M)\}$ be two $(H,M)$-orthogonal sets. Let $\cZ:=\inf(\cY^1,\cY^2)$ be the set of $Z_P:=\inf^P\{Y^1_P,Y^2_P,\}$ when $P$ runs $\cP(M)$. 
\begin{enumerate}\item The set $\cZ$ is a $(H,M)$-orthogonal set. 
\item  If $d(\cY^j)> 0$ for $j=1,2$ then   $d(\cZ)> 0$. In this case, the convex hull $\cS_M(\cZ)$ is the intersection of $\cS_M(\cY^1)$ and   $\cS_M(\cY^2)$.
\end{enumerate}
\end{lem}
\dem {\it 1.} Let $P$ and $P'$ two adjacent parabolic subgroups in $\cP(M)$ whose chambers in $a_M$  share the wall determined by the simple root $\al\in\De_P\cap(-\De_{P'})$.  
Let $\ga\in\De_P-\{\al\}$. By definition of orthogonal sets, for $j=1$ or $2$, one has $(\tilde{\om}_\ga^P,Y^j_P)=(\tilde{\om}_\ga^P,Y^j_{P'})$. By   Lemma \ref{copoids},  we have $\tilde{\om}_\ga^P=\tilde{\om}_{\ga'}^{P'}$. Hence we obtain 
 $(\tilde{\om}_\ga^P,Z_P)=(\tilde{\om}_{\ga'}^{P'},Z_{P'})$  and $(\tilde{\om}_{\ga'}^{P'},Z_{P'})=(\tilde{\om}_{\ga}^{P},Z_{P'})$. 
Since the scalar product on $a_0$ identifies $a_M$ to $a_M^*$,   one deduces that $Z_P-Z_{P'}$ is proportional to $\check{\al}$. 
\me

\no{\it 2.} Let $j\in\{1,2\}$ and $P\in\cP(M)$. By definition, we have $d(\cY^j)>0$ if and only if $\al(Y^j_P)> 0$ for all  $\al\in\De_P$.  By (\cite{ArDS} Corollary 2.2), this implies that $(\tilde{\om}_\al^P, Y^j_P)> 0$ for all  $\al\in\De_P$. Let $\al\in\De_P$. Writing $$Y_P^j=(\tilde{\om}_\al^P, Y^j_P) \al+\sum_{\be\in\De_P-\{\al\}} (\tilde{\om}_\be^P, Y^j_P) \be +X^j,$$ 
with $X^j\in a_H$, the condition $\al(Y^j_P)> 0$ is equivalent to 
 $$  \sum_{\be\in\De_P-\{\al\}}(\tilde{\om}_\ga^P,Y^j_P)[ -(\be, \al)]<( \tilde{\om}^P_\al,Y^j_P)(\al,\al).$$
 
 Since the real numbers $(\tilde{\om}_\be^P, Y^j_P)$ for $\be\in \De_P$ and $-(\be,\al)$ for $\al\neq \be$ in $\De_P$ are nonnegative, one deduces that 
 $$\sum_{\be\in\De_P-\{\al\}}(\tilde{\om}_\be^P, Z_P)[ -(\be, \al)]=\sum_{\be\in\De_P-\{\al\}}\inf\big( (\tilde{\om}_\be^P, Y^1_P), (\tilde{\om}_\be^P, Y^2_P)\big)[ -(\be, \al)] $$
 $$\leq \inf\big( \sum_{\be\in\De_P-\{\al\}}(\tilde{\om}_\be^P, Y^1_P)[-(\be, \al)], \sum_{\be\in\De_P-\{\al\}}(\tilde{\om}_\be^P, Y^2_P)[-(\be, \al)]\big)$$
 $$<\inf\big( (\tilde{\om}^P_\al, Y^1_P), (\tilde{\om}^P_\al, Y^2_P)\big)(\al,\al) =(\tilde{\om}^P_\al, Z_P)(\al,\al).$$
One deduces that $\al(Z_P)> 0$ for $\al\in\De_P$, thus
 $d(\cZ)> 0$.  \me

For the property of the convex hulls, it is enough to prove that, for all $P\in \cP(M)$, one has $a_P^+\cap \cS_M(\cY^1)\cap \cS_M(\cY^2)= a_P^+\cap \cS_M(\cZ)$. By (\cite{ArLT}, Lemma 3.1), one has $$  a_P^+\cap \cS_M(\cY^j)=\{ X\in a_P^+; \om_\ga^P(X-Y_P^j)\leq 0, \ga\in \De_P\}. $$ 
Since  $\tilde{\om}^P_\ga=c_\ga \om_\ga^P$ for $\ga\in\De_P$, where  $c_\ga$ is a positive real number,  the assertion follows easily.\qed

\subsection{ The truncated kernel}
We consider the regular representation $R$ of $G\times G$ on $L^2(G)$ defined by
$$ \big(R(y_1,y_2) \phi\big)(x)= \phi(y_1^{-1} x y_2),\quad \phi\in L^2(G), y_1, y_2\in G.$$
Consider $f\in C_c^\infty (G\times G)$ of the form $f(y_1, y_2)= f_1(y_1)f_2(y_2)$ with $f_j\in C_c^\infty(G)$. Then $$R(f):=\int_G\int_G f_1(y_1)f_2(y_2) R(y_1, y_2) dy_1 dy_2$$ is an integral operator with smooth kernel 
$$K_f(x,y)=\int_G f_1(xg) f_2(gy) dg= \int_G f_1(g) f_2(x^{-1} g y) dg.$$

In our case ($H$ is split), one has $A_H=A_G$, and the kernel $K_f$ is invariant by the diagonal $diag(A_H)$ of $A_H$. 
Since $H$ is not  compact, we introduce truncation to integrate this kernel on $diag(A_H)\bb(H\times H)$. \me

We fix a point $T$ in $a_{0,\F}$. If $P_0\in\cP(A_0)$, let $T_{P_0}$ be the unique translate by the Weyl group $W(H,A_0)$ of $T$ in the closure $\bar{a}_{P_0}^+$ of the positive Weyl chamber $a_{P_0}^+$. Then
$$\cY_T:=\{ T_{P_0}; P_0\in\cP(A_0)\}$$
is a   $(H,A_0)$-orthogonal set. We shall assume   that the number $$d(T):=\inf_{\al\in \De_{P_0}, \; P_0\in\cP(A_0)}\al(T_{P_0})$$  is suitable large. This means that the distance from $T$ to any of the root  hyperplanes in $a_0$ is large.  \ber\label{uxT} 
We denote by $u(\cdot,T)$ the characteristic function in $A_H\bb H$ of the set of points $x$ such that 
$$x=k_1 ak_2  \mbox{ with } a\in A_H\bb A_0, k_1,k_2\in K  \mbox{ and } h_{A_0}(a)\in \cS_{A_0}(\cY_T),$$
where $H=K A_0K  $ is the Cartan decomposition of $H$.\eer
We  consider $u(\cdot,T)$ as a $A_H$-invariant function on $H$. Thus, there is a compact set $\Om_T$ of $H$ such that if $u(x,T)\neq 0$ then $x\in A_H\Om_T$.  Let $\Om$ be a compact subset of $G$ containing the support  of $f_1$ and $f_2$. We consider $g\in G$ and $x_1, x_2\in H$ such that $f_1(g)f_2(x_1^{-1}gx_2)u(x_1,T) u(x_2,T)\neq 0$. Thus, there are $\om_1,\om_2$ in $\Om_T$ and $a_1, a_2$ in $A_H$ such that $x_1=\om_1 a_1$, $x_2=\om_2 a_2$ and  we have $g\in \Om$ and $x_1^{-1} gx_2=\om_1^{-1} g\om_2 a_1^{-1} a_2\in \Om$ since $A_H=A_G$. One deduces that $a_1^{-1} a_2$ lies a compact subset of $A_H$. Therefore the map $(g, x_1,x_2)\mapsto f_1(g)f_2(x_1^{-1}gx_2)u(x_1,T) u(x_2,T)$ is a compactly supported function on $G\times diag(A_H)\bb(H\times H)$.

\no Hence, we can define
$$ K^T(f):=\int_{diag(A_H)\bb H\times H} K_f(x_1, x_2) u(x_1,T) u(x_2, T) d\overline{(x_1,x_2)} .$$

\no By Fubini's Theorem, we have 

$$ K^T(f)=\int_G\int_{diag(A_H)\bb H\times H}  f_1(g) f_2(x_1^{-1} g x_2) u(x_1,T) u(x_2, T) d\overline{(x_1,x_2)} dg.$$\no We apply the Weyl integration formula (\ref{Weyl2}). Thus, we obtain
\beq\label{KT} K^T(f)=\sum_{M\in\cL(A_0)} c_M\sum_{S\in\cT_M}\sum_{x_m\in\kappa_S} c_{S,x_m}\int_{S_\si} K^T(x_m, \ga, f) d\ga,\eeq
where, for $S\in \cT_M$, $x_m\in\kappa_S$ and $\ga\in S_\si$, we have
$$K^T(x_m, \ga,f)=|\De_\si(x_m\ga)|_\F^{1/2}\int_{diag(A_M)\bb H\times H} \int_{diag(A_H)\bb H\times H} f_1(y_1^{-1}x_m\ga y_2) $$
$$\times f_2(x_1^{-1}y_1^{-1} x_m\ga y_2 x_2) u(x_1,T) u(x_2, T) d\overline{(x_1,x_2)} d\overline{(y_1,y_2)}.$$
We recall that each $x_m$ in $\kappa_S$ and $\ga$ in $S_\si$  commute with $A_M$  for $S\in\cT_M$.\\
We first replace  $(x_1,x_2)$ by $(y_1 x_1, y_2x_2)$ in the integral over  $\overline{(x_1,x_2)}$. The resulting integral over $diag(A_H)\bb H\times H$ can be expressed as a double integral over $a\in A_H\bb A_M$ and $(x_1,x_2)\in diag(A_M)\bb H\times H$ which depends on  $\overline{(y_1,y_2)}\in diag(A_M)\bb H\times H$. Since $A_M$ commutes with $x_m\in\kappa_S$ and $\ga\in S_\si$, we obtain


\ber\label{KTga} $$K^T(x_m, \ga,f)=|\De_\si(x_m\ga)|_\F^{1/2}\int_{diag(A_M)\bb H\times H} \int_{diag(A_M)\bb H\times H} f_1(y_1^{-1}x_m\ga y_2) $$
$$\times f_2(x_1^{-1} x_m\ga   x_2) u_M(x_1,y_1,x_2,y_2, T)   d\overline{(x_1,x_2)} d\overline{(y_1,y_2)}$$ where 
$$u_M(x_1,y_1,x_2,y_2, T) =\int_{A_H\bb A_M} u(y_1^{-1}ax_1,T) u(y_2^{-1}ax_2, T) da.$$\eer

Our goal is to prove that $K^T(f)$ is asymptotic to an expression $J^T(f)$ where $ J^T(f)$ is obtained in a similar way to $K^T(f)$ where we replace the weight function $u_M(x_1,y_1,x_2,y_2, T) $ by another weight function $v_M(x_1,y_1,x_2,y_2, T) $ defined as follows.

We fix 
$M\in\cL(A_0)$ and   $P\in \cP(M)$. Let  $P_0\in\cP(A_0)$ be contained in $P$. We denote by  $T_P$ the projection of $T_{P_0}$ on $a_M$ according to the decomposition $a_0=a_M\oplus a_0^M$.   By (\ref{restY}), the set $\cY_M(T):=\{T_P; P\in\cP(M)\}$ is a $(H,M)$-orthogonal set independent of the choices of $P_0$.  Moreover, by (\cite{ArLT} (3.2)), we have   $d(\cY_M(T))\geq d(T)>0$. Thus, $\cY_M(T)$ is   positive.\me

For $x, y$ in $H$, we set 
$$Y_P(x,y,T):= T_P+h_P(y)-h_{\overline{P}}(x).$$

By (\cite{ArLT}, page 30), the set $\cY_M(x,y,T):=\{ Y_P(x,y,T); P\in\cP(M)\}$ is a $(H,M)$-orthogonal set, which is positive when $d(T)$ is sufficiently large relative to $x$ and $y$. \me

For $x_1, x_2, y_1$ and $y_2$ in $H$, we set
  \beq\label{ZP} Z_P(x_1,y_1,x_2,y_2,T):=\mbox{inf}^P( Y_P(x_1,y_1,T),  Y_P(x_2,y_2,T))\eeq
where $\mbox{inf}^P$ is defined in (\ref{infP}) and 
\beq\label{YxyT} \cY_M(x_1,y_1,x_2,y_2,T):= \{ Z_P(x_1,y_1,x_2,y_2,T); P\in\cP(M)\}.\eeq

By Lemma \ref{infP}, the set $\cY_M(x_1,y_1,x_2,y_2,T)$ is a $(H,M)$-orthogonal set. Moreover,
when   $d(T)$ is large   relative to $x_i, y_i$, for $i=1,2$, one has $d(\cY_M(x_1,y_1,x_2,y_2,T))> 0$, hence  this set is positive. 
We   define the weight function $v_M$ by
\beq\label{vM} v_M(x_1,y_1,x_2,y_2, T) := \int_{A_H\bb A_M}\si_M(h_M(a), \cY_M(x_1, y_1, x_2, y_2, T)) da\eeq
where $\si_M$ is defined in (\ref{SigmaM}).

We set
\beq\label{JT} J^T(f):=\sum_{M\in\cL(A_0)} c_M\sum_{S\in\cT_M}\sum_{x_m\in\kappa_S} c_{S,x_m}\int_{S_\si} J^T(x_m,\ga,f) d\ga,\eeq
where
\ber\label{JTga}$$J^T(x_m, \ga,f)=|\De_\si(x_m\ga)|_\F^{1/2}\int_{diag(A_M)\bb H\times H} \int_{diag(A_M)\bb H\times H} f_1(y_1^{-1}x_m\ga y_2) $$
$$\times f_2(x_1^{-1} x_m\ga   x_2) v_M(x_1,y_1,x_2,y_2, T)   d\overline{(x_1,x_2)} d\overline{(y_1,y_2)}.$$\eer

Our main result is the following. We will prove it in section \ref{preuve}.

\begin{theo}\label{resultatprincipal} Let $\de>0$. Then, there are positive numbers $C$ and $\ep$ such that for all $T$ with $d(T)\geq \de\Vert T\Vert$, one has
\beq\label{KT-JT}\vert K^T(f)-J^T(f)\vert\leq C e^{-\ep\Vert T\Vert}.\eeq
\end{theo}
\subsection{Preliminaries to estimates} 
We fix a norm  $\Vert\cdot \Vert$ on $G$  as in (\ref{norme}).
Let  $\F'$ be a finite extension of $\F$. We set $\sG':=\sG\times_\F\F'$ and $G':=\sG'(\F')$.  One can extend the absolute value $\vert\cdot \vert_\F$ to $\F'$, and the   norm $\Vert\cdot \Vert$ to  $G'$.  For $x,y$ in $G'$, we set 
$$\Vert (x,y)\Vert:=\Vert x\Vert \Vert y\Vert.$$

To obtain our estimates, we will use notation of (\ref{comparfg}) and (\ref{equival}). Since the norm takes values greater or equal to $1$, we will  freely apply the properties (\ref{compar>1}).
\begin{lem}\label{estimss'} Let 
$S$ be a maximal torus of $H$ and  let $M$ be the centralizer of $A_S$ in $H$. We fix $x_m\in G\cap\sM\sS_\si=\tilde{M}\cap \sM\sS_\si$.  Then, one has
\beq\label{ss'}  \inf_{s\in S} \Vert (sx_m^{-1}x_1, sx_2)\Vert\preccurlyeq \inf_{s'\in \sS(\F')} \Vert (s'x_m^{-1}x_1, s'x_2)\Vert,\quad x_1, x_2\in H.\eeq
\end{lem}

\dem Since $H^1 A_H$ is of finite index in $H$, using (\ref{compact}) we may assume that $x_1,x_2$ belong to $H^1 A_H$. Since $A_G=A_H$, using the invariance of the property by the left action of $diag(A_H)$ on $(x_1,x_2)$, it is enough to prove the result for  $x_1\in H^1$ and $x_2=a_2 y_2$ with $a_2\in A_H$ and $y_2\in H^1$.

To establish (\ref{ss'}), we first assume that $A_S=A_H$ which implies that  the quotient $A_H\bb S$ is compact. By (\ref{compact}), there is a positive constant $C$ such that
$$ \inf_{s\in S} \Vert (sx_m^{-1}x_1, sx_2)\Vert\leq C \inf_{a\in A_H} \Vert (ax_m^{-1}x_1,ax_2)\Vert.$$
We deduce from (\ref{prod}) that
$$\Vert (ax_m^{-1}x_1,ax_2)\Vert\leq\Vert x_m^{-1}\Vert \Vert a\Vert^2 \Vert a_2\Vert \Vert x_1\Vert \Vert y_2\Vert.$$

\no Taking  the lower bound in $a\in A_H$,   there is a positive constant $C_1$ such that
\beq\label{estim1} \inf_{s\in S} \Vert (sx_m^{-1}x_1, sx_2)\Vert\leq C_1  \Vert x_1\Vert\Vert a_2\Vert  \Vert y_2\Vert.\eeq
We now use the following  Lemma of \cite{ArLT} (Lemma   4.1):
\ber\label{ar41} If $S_0$ is a maximal  torus of $H$ with $A_H\bb S$ compact, then there exists an element $s_0\in S_0$
such that $$\Vert y\Vert \preccurlyeq\Vert y^{-1} s_0y\Vert, \quad y\in H^1.$$\eer
We apply  this Lemma to $S_0=S$. Since $\sS(\F')$ commutes with $s_0$, using  the property (\ref{prod}) of the norm, one deduces
\beq\label{yF'} \Vert y_2\Vert\preccurlyeq \Vert s' y_2\Vert^2\Vert s_0\Vert,\quad y_2\in H^1, s'\in \sS(\F').\eeq
 On the other hand $S_1:=x_m Sx_m^{-1}$ is a  maximal torus of $H$ which satisfies $A_{S_1}=A_H$ since $x_m\in G\cap\sM\sS_\si$. Applying (\ref{ar41}) to $S_0=S_1$, there exists $s_1\in S$ such that 
\beq\label{x41} \Vert x_1\Vert\preccurlyeq \Vert x_1^{-1} x_m s_1 x_m^{-1} x_1\Vert,\quad x_1\in H^1.\eeq
The same argument as above leads to 
\beq\label{xF'} \Vert x_1\Vert\preccurlyeq \Vert s' x_m^{-1}x_1\Vert^2\Vert s_1\Vert, \quad x_1\in H^1, s'\in \sS(\F').\eeq\\
Then, by (\ref{estim1}), (\ref{yF'}) and (\ref{xF'}), and applying the properties  (\ref{compar>1}), we deduce that  
\beq\label{estim2}  \inf_{s\in S} \Vert (sx_m^{-1}x_1, sa_2 y_2)\Vert\preccurlyeq \Vert s' x_m^{-1} x_1\Vert \Vert s' y_2\Vert\Vert a_2\Vert,\quad s'\in \sS(\F'), x_1, y_2\in H^1, a_2\in A_H.\eeq\\
To obtain our result, we have to prove that 
\beq\label{finestim}\Vert s' x_m^{-1} x_1\Vert \Vert s' y_2\Vert\Vert a_2\Vert\preccurlyeq \Vert (s' x_m^{-1}x_1, s'a_2y_2)\Vert,\quad s'\in \sS(\F'), x_1,y_2\in H^1, a_2\in A_H.\eeq\\
We can write  $\sS=\sT \sA_H$ where  $\sT$  is a  maximal torus of the derived group $\sH_{der}$ of $\sH$. We set  $T':=\sT(\F')$ and  $A'_H:=\sA_H(\F')$. Then $T'$ is contained in $ H'^1$. Moreover, the intersection  of $\sT$ and $\sA_H$ is finite. Hence, one has the exact sequence 
$$1\to \sT\cap \sA_H\to \sT\times \sA_H\to \sS\to 1.$$
 Going to $\F'$-points, the long exact sequence in cohomology implies that   $T'A'_H$ is of finite index in $\sS(\F')$. By (\ref{compact}), it is enough to prove (\ref{finestim}) for  $s'=t'a'\in \sS(\F')$ with $t'\in T'$ and $a'\in A'_H$. By (\ref{G1}), if $x_1\in H^1$ then    $x_1\in   H'^1\subset G'^1$ and  $x_m^{-1}x_1x_m\in G'^1$. Since $H$ is split, we have   $A'_H=A'_G$. Then   (\ref{g1a}) gives 
  $$\Vert a't' x_m^{-1}x_1\Vert\approx \Vert a't' x_m^{-1}x_1 x_m\Vert\approx \Vert a'\Vert \Vert t' x_m^{-1} x_1x_m\Vert, \quad a'\in A'_H, t'\in T', x_1\in H^1, $$ \me
and 
$$\Vert a't'y_2\Vert\approx\Vert a'\Vert \Vert t' y_2\Vert \quad a'\in A'_H, t'\in T', y_2\in H^1.$$\\
Applying (\ref{compar>1}), we deduce that
\beq\label{Ia't'}\Vert t'a' x_m^{-1} x_1\Vert \Vert a't'y_2\Vert\Vert a_2\Vert\approx \Vert a_2\Vert \Vert a'\Vert^2\Vert t' x_m^{-1} x_1x_m\Vert \Vert t'y_2\Vert \approx  \Vert a_2\Vert \Vert a'\Vert \Vert t' x_m^{-1} x_1x_m\Vert \Vert t'y_2\Vert,\eeq
for $ t'\in T', a'\in A'_H, x_1,y_2\in H^1, a_2\in A_H .$

\no Let us prove that \beq\label{aa_2} \Vert a' \Vert \Vert a' a_2\Vert\approx \Vert a'\Vert \Vert a_2\Vert,\quad a'\in A'_H, a_2\in A_H.\eeq

\no We have $\Vert a' a_2\Vert\leq \Vert a'\Vert \Vert a_2\Vert$ by (\ref{prod}). Then   $\Vert a'\Vert \Vert a' a_2\Vert\leq (\Vert a'\Vert \Vert a_2\Vert)^2$ since $1\leq  \Vert a_2\Vert$. As $\Vert a'\Vert=\Vert a' a_2 a_2^{-1}\Vert\leq \Vert a' a_2\Vert \Vert a_2\Vert$, we have $\Vert a'\Vert \Vert a_2\Vert\leq (\Vert a'a_2\Vert \Vert a_2\Vert)^2$ and (\ref{aa_2}) follows. Applying  (\ref{aa_2}) in (\ref{Ia't'}), we deduce that
\ber\label{majI'} $$\Vert t'a' x_m^{-1} x_1\Vert \Vert a't'y_2\Vert\Vert a_2\Vert\preccurlyeq \Vert a'\Vert \Vert t' x_m^{-1} x_1 x_m\Vert \Vert a' a_2\Vert \Vert t'y_2\Vert,$$
for $ t'\in T', a'\in A'_H, x_1,y_2\in H^1, a_2\in A_H.$\eer
Since $x_m^{-1} H^1x_m \subset G'^1$ and $A'_H=A'_G$, we obtain from  (\ref{g1a})   $$\Vert a'\Vert\Vert t' x_m^{-1} x_1 x_m\Vert\approx \Vert a'  t' x_m^{-1} x_1 x_m\Vert\approx \Vert a'  t' x_m^{-1} x_1 \Vert, \quad a'\in A'_H, t'\in T', x_1\in H^1, $$ and $$\Vert a' a_2\Vert \Vert t'y_2\Vert\approx \Vert a' a_2 t'y_2\Vert, \quad  a'\in A'_H, t'\in T', a_2\in A_H, y_2\in H^1.$$\\
  Applying this in (\ref{majI'}) and using (\ref{compar>1}), we deduce that
\ber$$\Vert t'a' x_m^{-1} x_1\Vert \Vert a't'y_2\Vert\Vert a_2\Vert\preccurlyeq \Vert a'  t' x_m^{-1} x_1 \Vert \Vert a' t'a_2 y_2\Vert$$
for  $a'\in A'_H$, $t'\in T'$, $x_1,y_2\in H^1$.\eer
 Then, the property (\ref{finestim}) follows. This finishes the proof of the Lemma when $A_H\bb S$ is compact.\me

We now prove (\ref{ss'}) for any maximal torus $S$ of $H$. Let $A_S$ be the maximal split torus of $S$ and $M$ be the centralizer of $A_S$ in $H$. Thus we have $A_M=A_S$ and $A_M\bb S$ is  compact. Let $P=MN_P\in\cP(M)$ and let $K$ be a compact subgroup of $H$ such that $H=PK$. Each $x\in H$ can be written $x=m_P(x) n_P(x) k(x)$ with $m_P(x)\in M, n_P(x)\in N_P$ and $k(x)\in K$. Then, there is a positive constant $C$ such that 
\beq\label{infmP}  \inf_{s\in S} \Vert (sx_m^{-1}x_1, sx_2)\Vert\leq  C\inf_{s\in S} \big(\Vert sx_m^{-1}m_P(x_1)\Vert \Vert sm_P(x_2)\Vert\big)\Vert n_P(x_1)\Vert \Vert n_P(x_2)\Vert,\eeq
for $x_1, x_2\in H$.
By assumption on $x_m$, there are $h_m\in\sM$ and $a_m\in \sS_\si$ such that $x_m=h_m a_m\in \tilde{M}$. Hence, we can applied the first part of the proof to $(M,S)$ instead of $(H,S)$.  Therefore, we  obtain 
$$  \inf_{s\in S} \Vert (sx_m^{-1}x_1, sx_2)\Vert\preccurlyeq   \inf_{s'\in \sS(\F')} \big(\Vert s'x_m^{-1}m_P(x_1)\Vert \Vert s'm_P(x_2)\Vert\ \big)\Vert n_P(x_1)\Vert \Vert n_P(x_2)\Vert, \quad x_1,x_2\in H.$$
To compare the right hand-side of this inequality to those of (\ref{ss'}), we will use  the Iwasawa decomposition (\ref{Iwasawa}) of $H'$. Let $K'$  be a compact subgroup of $H'$ such that  $H'= \sP(\F')K'=\sM(\F')\sN_P(\F') K'$. According to (\ref{mP}), each  $y$ in $ H'$ can be written $y=m'_P(y)n'_P(y)k'$ with   $m'_P(y)\in \sM(\F')$, $n'_P(y)\in \sN_P(\F')$ and $k'\in K'$. Then for $x\in H$ and $z\in \sM(\F')$, we have $zx=zm_P(x)n_P(x) k=m'_P(zx) n'_P(zx)k'$ with $k\in$ and $k'\in K'$. Hence,  since $K$ and $K'$ are compact subsets, there is a positive constant $C'$ such that 
$$\Vert n'_P(zx)^{-1} m'_P(zx)^{-1}zm_P(x)n_P(x)\Vert  \leq C', \quad z\in \sM(\F'), x\in H.$$
 Since $zm_p(x)\in \sM(\F')$ for $z\in\sM(\F')$ and $x\in H$, we deduce from (\ref{mn}) that there is a positive constant $C_1$ such that for $x\in H$ and $z\in \sM(\F')$, one has
$$\Vert n'_P(zx)^{-1} m'_P(zx)^{-1}zm_P(x) n'_P(zx)\Vert  \leq C_1\quad\mbox{and } \quad \Vert n'_P(zx)^{-1}n_P(x)\Vert  \leq C_1.$$ By (\ref{prod}), we obtain
$$ \Vert zm_P(x) \Vert  \leq C_1\Vert m'_P(zx)\Vert  \Vert n'_P(zx)\Vert^2\quad\mbox{and }\quad\Vert n_P(x)\Vert\leq C_1\Vert n'_P(zx)\Vert .$$
Using (\ref{mn}) again, it follows that
$$\quad \Vert zm_P(x) \Vert \preccurlyeq  \Vert zx\Vert,\quad\mbox{and }\Vert n_P(x)\vert\preccurlyeq  \Vert zx\Vert,\quad z\in \sM(\F'), x\in H,$$
hence by (\ref{compar>1})

\beq\label{zmP}\Vert zm_P(x) \Vert\Vert n_P(x)\vert\preccurlyeq  \Vert zx\Vert,\quad z\in \sM(\F'), x\in H.\eeq
We deduce that 
\beq\label{x2}\Vert s'm_P(x_2) \Vert\Vert n_P(x_2)\vert\preccurlyeq  \Vert s'x_2\Vert,\quad s'\in \sS(\F'), x_2\in H.\eeq
Since  $x_m=h_m a_m$ with $h_m\in\sM$ and $a_m\in \sS_\si$, one has $x_m s'x_m^{-1}\in\sM\cap H'=\sM(\F')$ for $s'\in \sS(\F')$. Therefore, we deduce from (\ref{zmP}) that 
\beq\label{x1}\Vert x_ms'x_m^{-1}m_P(x_1) \Vert\Vert n_P(x_1)\vert\preccurlyeq  \Vert x_ms'x_m^{-1}x_1\Vert,\quad s'\in \sS(\F'), x_1\in H.\eeq
Since $\Vert s'x_m^{-1}m_P(x_1)\Vert\leq \Vert x_m^{-1}\Vert \Vert x_ms'x_m^{-1}m_P(x_1)\Vert$ and  $ \Vert x_ms'x_m^{-1}x_1\Vert \leq  \Vert x_m\Vert\vert s'x_m^{-1}x_1\Vert$, we deduce the estimate (\ref{ss'}) from (\ref{infmP}), (\ref{x2}) and (\ref{x1}). 
This finishes the proof of the Lemma.\qed\me

The following Lemma is the analogue of Lemma 4.2 of (\cite{ArLT}).
\begin{lem}\label{inf-delta}  Let $S$ be a maximal torus of $H$ and let  $x_m\in \kappa_S$. Then, there is a positive integer $k$ with the property  that, for any given  compact subset $\Om$ of $G$, there exists a positive constant $C_\Om$  such that, for all $\ga\in S_\si$ with $x_m\ga\in G^{\si-reg}$, and all $x_1, x_2$ in $H$ satisfying $x_1^{-1} x_m\ga x_2\in \Om$, one has 
$$\inf_{s\in S} \Vert (sx_m^{-1}x_1, sx_2)\Vert\leq C_\Om \vert\De_\si(x_m\ga)\vert_\F^{-k}.$$
\end{lem}
\dem 
Let $\F'$ be a finite extension   of   $\E$ such that $\tilde{\sS}$  splits over $\F'$. Recall that we can write $x_m=h_m a_m$ with $h_m\in\sH$ and $a_m\in \sS_\si$. Thus we may and will choose $\F'$ such  that $h_m\in\sH(\F')$ and $a_m\in\sS_\si(\F')$. For convenience of lecture, if $\sJ$ is an algebraic variety defined over $\F$, we set $J':=\sJ(\F')$.

By the previous Lemma \ref{estimss'}, it is enough to prove the existence of  a positive integer $k$ satisfying the property that for any compact subset $\Om'$ of $G'^{\si-reg}$, there exists  $C_{\Om'}>0$ such that   
\beq\label{inf-delta'}  \inf_{s'\in S'} \big( \Vert s'x_m^{-1}x_1\Vert \Vert  s'x_2\Vert\big)\leq C_{\Om'} \vert\De_\si(x_m\ga)\vert_\F^{-k}\eeq
 for all $x_1, x_2\in H'$ and  $\ga\in S_\si$ satisfying  $x_m\ga\in G^{\si-reg}$ and $x_1^{-1} x_m\ga x_2\in \Om'$.\me

Let $B'=S' N'$ be a Borel subgroup of $H'$ containing $S'$ and $K'$ be  a compact subgroup of $H'$ such that $H'=S' N' K'=N'S'K'$. We can also write $H'=(h_mS'h_m^{-1})(h_m N' h_m^{-1}) (h_m K' h_m^{-1})$. By (\ref{compact}), one can reduce the proof of the statement for $x_1\in (h_mS'h_m^{-1})(h_m N' h_m^{-1}) $ and $x_2\in S' N'$. 


Let  $x_1=h_m s_1 n_1 h_m^{-1}$ and $x_2=s_1 s_2 n_2$ with $s_1, s_2\in S'$ and $n_1, n_2\in N'$. Since $x_m= h_m a_m$, we have $x_ms_1x_m^{-1}=h_ms_1h_m^{-1}$, hence  for $s'\in S'$, we have  $ s'x_m^{-1}x_1=s'x_m^{-1}x_ms_1x_m^{-1}h_m n_1 h_m^{-1}=s's_1x_m^{-1} h_m n_1 h_m^{-1}$. We obtain 

$$\inf_{s'\in S'} \big( \Vert s'x_m^{-1}x_1\Vert \Vert  s'x_2\Vert\big)=\inf_{s'\in S'} \big( \Vert s'x_m^{-1}h_m n_1 h_m^{-1}\Vert \Vert  s's_2\Vert\big).$$\\
Notice that $x_1^{-1} x_m\ga x_2=h_m n_1^{-1}h_m^{-1}x_ms_1^{-1}x_m^{-1}x_m\ga s_1s_2n_2=h_m n_1^{-1}h_m^{-1}x_m\ga s_2n_2$. 

Therefore, we are reduced to   prove (\ref{inf-delta'}) for  $x_1= h_m n_1 h_m^{-1}$ with $n_1\in N'$,    $x_2=s_2n_2$ with $n_2\in N'$, $s_2\in S'$ and $\ga\in S_\si $ such that $x_m \ga$ is $\si$-regular and $x_1^{-1} x_m\ga x_2\in \Om'$.  By the properties of the norm,  there is some positive constant $C'$ such that 
\beq\label{reducpreuve} \inf_{s'\in S'} \big( \Vert s'x_m^{-1}x_1\Vert \Vert  s'x_2\Vert\big)\leq C'\Vert n_1\Vert \Vert s_2\Vert \Vert n_2\Vert, \quad x_1=h_m n_1 h_m^{-1}, x_2=s_2 n_2.\eeq
We want to estimate $\Vert n_1\Vert \Vert s_2\Vert \Vert n_2\Vert$  when $x_1= h_mn_1h_m^{-1} $ and $x_2= s_2 n_2$ satisfy $x_1^{-1} x_m \ga x_2\in \Om'$.\
 For this, we use the isomorphism $\Psi$  from $G'$ to $H'\times H'$ defined in (\ref{isoProduit}). If $x\in H'$ then   $\Psi(x)=(x,  x)$   and if  $y\in G$ satisfies $y^{-1}=\si(y)$ then $\Psi(y)=(y, y^{-1})$. We set $(y_1, y_2):=\Psi(x_1^{-1} x_m \ga x_2)$. Then, we have
$$y_1= h_m n_1^{-1} a_m \ga n_2 s_2= h_m (n_1^{-1} a_m\ga n_2 (a_m\ga)^{-1}) (a_m\ga s_2), $$ and $$ y_2=h_m n_1^{-1} a_m^{-1} \ga^{-1} n_2s_2=h_m \big( n_1^{-1} a_m^{-1} \ga^{-1}n_2 \ga a_m\big) (a_m\ga)^{-1} s_2.$$
Since $H'=N'S'K'$, the condition $x_1^{-1} x_m\ga x_2\in \Om'$ implies that there exist two  compact subsets $\Om_N\subset N'$ and $\Om_S \subset S'$ depending only on $\Om'$ such that 
$$n_1^{-1} a_m\ga n_2 (a_m\ga)^{-1}\in\Om_N,\quad\textrm{and}\quad  n_1^{-1} a_m^{-1} \ga^{-1}n_2 \ga a_m\in\Om_N,$$
$$a_m\ga s_2\in \Om_S,\quad\textrm{and}\quad(a_m\ga)^{-1} s_2\in \Om_S.$$
We deduce from the second property that $s_2$ and $\ga$ must lie in compact subsets of $S'$. 
We set $$\nu_1(\ga, n_1,n_2):=n_1^{-1} a_m\ga n_2 (a_m\ga)^{-1}\quad\mbox{ and }\quad\nu_2(\ga, n_1, n_2):= n_1^{-1} (a_m\ga)^{-1} n_2 a_m\ga.$$ We consider the map $\psi$ from $N'\times N'$ into itself defined by $\psi(n_1,n_2)=(\nu_1,\nu_2)$. Recall that $\Phi(S',\g h')$ denotes the set of roots of $S'$ in the Lie algebra $\g h'$ of $H'$ (cf.  \ref{roots}).  Let $\g n'$ be the Lie algebra of $N'$. For $\al\in \Phi(S',\g h')$, we denote by $X_\al\in \g n'$  the root vector in $\g h'$  corresponding to $\al$. Then $(a_m\ga)$ acts on $X_\al $ by $a_\al:=(a_m\ga)^\al$.  The differential $d_{(n_1,n_2)}\psi$ of $\psi$ at   $(n_1,n_2)\in N'\times N'$ is given by 
$d_{(n_1,n_2)}\psi (X_1, X_2)=(\Ad(a_m\ga n_2^{-1}(a_m\ga)^{-1})Y_1, \Ad((a_m\ga)^{-1} n_2^{-1}a_m\ga)Y_2)$ where $$Y_1=- \Ad(n_1)X_1+\Ad(a_m\ga)\Ad(n_2) X_2 $$ and $$Y_2=- \Ad(n_1)X_1+\Ad(a_m\ga)^{-1}\Ad(n_2) X_2.$$ 
The map $(X_1,X_2)\mapsto (Y_1,Y_2)$ is  the composition of the map $(X_1,X_2)\mapsto (\Ad(n_1) X_1, \Ad(n_2)X_2)$, whose determinant is equal to $1$, with  $d_e\psi$ where $e$ is the neutral point of $N'\times N'$. We deduce that 
the  jacobian of $\psi$ at $(n_1,n_2)$ is independent of $(n_1,n_2)$. At the neutral point $e\in N'\times N'$, we have 
 $ d_e\psi(X_\al, 0)=(-X_\al, -X_\al)$ and $d_e\psi(0, X_\al)=( a_\al  X_\al, a_{-\al}  X_\al)$. Hence, the jacobian of $\psi$ is equal to 
$$\vert \prod_{\al\in \Phi(S',\g h')} a_\al(1-a_{-2\al})\vert_{\F'}= \vert \mbox{det}(\Ad(a_m\ga))_{\g h'/\g s'}\vert_{\F'}\vert \mbox{det}(1-\Ad(a_m\ga)^{-2})_{\g h'/\g s'}\vert_{\F'}=\vert D_{H'}((a_m\ga)^{-2})\vert_{\F'}.$$
Recall that $x_m\ga$ is assumed to be $\si$-regular. Thus, by (\ref{deltas}), one has $\De_\si(x_m\ga)=D_{H'}(a_m^{-2}\ga^{-2})\neq 0$ . Then,  arguing as in (\cite{HC2Thm}  proof of Lemma 10 and Lemma 11), we deduce that  the map $\psi$  is an $\F'$-rational isomorphism of $\sN\times \sN$ to itself whose inverse  $(\nu_1,\nu_2)\mapsto (n_1, n_2):=(n_1(\ga, \nu_1,\nu_2), n_2(\ga, \nu_1, \nu_2))$ is rational.  Moreover, there is a positive integer $k$ such that the map 
$$(y,\nu_1,\nu_2)\mapsto D_{\sH}(y)^k (n_1(y, \nu_1,\nu_2),n_2(y, \nu_1,\nu_2))$$
is defined by an $\F'$-rational morphism between the algebraic varieties $\sS\times \sN\times\sN$ and $\sN\times \sN$. Since  $\nu_1,\nu_2
$ and $\ga$ lie in  compact subsets depending only on $\Om'$, one deduces that there exists a constant $C_{\Om'}>0$ such that 
$$\Vert (n_1(\ga,\nu_1,\nu_2), n_2(\ga,\nu_1,\nu_2))\Vert\leq C_{\Om'} \vert D_{H'}(a_m^{-2}\ga^{-2})\vert_{\F'}^{-k}= C_{\Om'} \vert\De_\si(x_m\ga)\vert_\F^{-k}.$$
The Lemma follows from  (\ref{reducpreuve}) and the fact that $s_2$ lies in a compact set. \qed
\subsection{Proof of Theorem \ref{resultatprincipal}}\label{preuve}
Our goal is to prove that $\vert K^T( f)-J^T(f)\vert $  is bounded by a function which approaches $0$ as $T$ approaches infinity. By definition, $K^T(f)$ and $J^T(f)$ are finite linear combinations  of $\int_{S_\si} K^T(x_m, \ga,f) d\ga$ and $\int_{S_\si} J^T(x_m, \ga,f) d\ga$ respectively, where $M\in\cL(A_0)$, $S$ is a maximal torus of $M$ satisfying  $A_S=A_M$ and $x_m\in \kappa_S$ (cf.  (\ref{KT}) and (\ref{JT})). \me

 We fix $M\in\cL(A_0)$ and a maximal  torus $S$ of $M$ such that $A_S=A_M$. Let $x_m\in\kappa_S$. To obtain our result, it is enough to establish the estimate (\ref{KT-JT}) for  $\int_{S_\si} \vert K^T(x_m,\ga,f)-J^T(x_m, \ga,T)\vert d\ga$. This will be done in the Corollary \ref{KTx-JTx} below.\me
 
For $\ep>0$, we define
\beq\label{SsepT} S_\si(\ep,T):=\{ \ga\in S_\si; 0<\vert \De_\si(x_m\ga)\vert_\F \leq e^{-\ep \Vert T\Vert }\}.\eeq
\begin{lem}\label{intDe}
\begin{enumerate}\item There exists $\ep_0>0$ such that the map $\ga\mapsto \vert\De_\si(x_m\ga)\vert_{\F}^{-\ep_0}$ is locally integrable on $S_\si$. 
\item Let $\ep>0$. Let $B$ be a bounded subset of $S_\si$ and  $p$ be a nonnegative integer. Then, there is a positive constant $C_{B,p}$ depending on $B$ and $p$, such that
$$\int_{S_\si(\ep,T)\cap B} |\log|\De_\si(x_m\ga)|_\F^{-1}|^p d\ga\leq C_{B,p} e^{-\frac{\ep\ep_0\Vert T\Vert}{2}}.$$
\end{enumerate}
\end{lem}
\dem {\it 1.}  The proof follows those of the group case, we use the similar statement on Lie algebras and the exponential map. We denote by $\g s$ the Lie algebra of $S$. For $X\in \g s$, we set $\eta(X)=|det(\ad X_{| \g h/\g s})|_\F$.  By (\cite{HCvD} Lemma 44), there exists $\ep_0>0$ such that $X\mapsto \eta(X)^{-2\ep_0}$ is locally integrable on $\g s$. 
To obtain the result, it is sufficient to prove that 
\ber\label{intconv}for each $\ga_0\in S_\si$, there exists a compact neighborhood $U_0$ of $1$ such that the integral $\displaystyle \int_{U_0} \vert \De_\si(x_m\ga_0\ga)\vert_\F^{-\ep_0} d\ga$ converges.\eer\\
 If $x_m\ga_0$ is $\si$-regular, there is a compact neighborhood $U_0$ of $1$ in $S_\si$ such that $\vert \De_\si(x_m\ga_0\ga)\vert_\F=\vert \De_\si(x_m\ga_0)\vert_\F\neq 0$  for all $\ga\in U_0$. Hence (\ref{intconv}) is clear.
 
We assume that $x_m\ga_0$ is not $\si$-regular.  We choose an extension $\F'$ of $\E$ such that $\tilde{\sS}$ splits over $\F'$ and $\sp(x_m)\in \tilde{\sS}_\si(\F')$. We use notation of (\ref{roots}). Let $\Phi_0$ be the set of root $\al$ in $\Phi(S'_\si,\g g')$ such that $\sp(x_m\ga_0)^\al=1$. We set
$$\nu(\ga)=\prod_{\al\in\Phi(S'_\si,\g g')-\Phi_0}\vert 1-\sp(x_m\ga_0)^{\al}\ga^{-2\al}\vert^2_{\F'}.$$
We have $\De_\si(x_m\ga_0\ga)=D_{G'}(\sp(x_m\ga_0)\ga^{-2})=det(1-\Ad \sp(x_m\ga_0)\ga^{-2})_{|\g g/\tilde{\g s}}$ and each root of  $\Phi(S'_\si,\g g')$ has multiplicity $2$. Hence, we obtain 
$$\vert \De_\si(x_m\ga_0\ga)\vert_{\F'}=\nu(\ga)\prod_{\al\in\Phi_0}\vert 1-\ga^{-2\al}\vert^2_{\F'}.$$
We choose a  compact neighborhood $W$ of $1$ in $S_\si$ such that  $\nu(\ga)=\nu(1)\neq 0$ for $\ga\in W$. Let $\displaystyle \be=\sup_{\ga\in W} \prod_{\al\in\Phi(S'_\si,\g g')-\Phi_0}\vert 1-\ga^{-2\al}\vert^2_{\F'}.$ Then, for $\ga\in W$, we have

$$\be \vert \De_\si(x_m\ga_0\ga)\vert_{\F'}=\be\nu(1) \prod_{\al\in\Phi_0}\vert 1-\ga^{-2\al}\vert^2_{\F'}\geq \nu(1) \vert\De_\si(\ga)\vert_{\F'}.$$
Consider the exponential map, there exist two open neighborhoods $\om$ and $U$ of  $0$ and $1$ in $\g s$ and $S_\si$ respectively, such that the map  $X\mapsto \exp(\tau X)$ is well-defined on $\om$ and is an isomorphism and an homeomorphism onto $U$. For $X\in\om$, we have
$$\frac{\vert\De_\si(\exp(\tau X)\vert^{1/2}_{\F'}}{\eta(X)}=\prod_{\al\in\Phi(S'_\si,\g g')}\frac{\vert 1-e^{2\tau\al(X)}\vert_{\F'}}{\vert\al(X)\vert_{\F'}}.$$

We can choose a compact neighborhood  $\om_0\subset \om$ of $0$ in $\g s$ such that  this product is a positive constant $c $ and $U_0:=\exp(\tau \om_0)$ is contained in $W$. We deduce that
$$\int_{U_0}\vert \De_\si(x_m\ga_0\ga)\vert^{-\ep_0}_{\F}d\ga\leq \Big(\frac{\nu(1)}{\be}\Big)^{-\ep_0}\int_{U_0} \vert\De_\si( \ga)\vert^{-\ep_0}_{\F}d\ga= \Big(\frac{\nu(1)}{\be}\Big)^{-\ep_0}c\int_{\om_0}\eta(X)^{-2\ep_0} dX.$$ 
The right hand side of this inequality is finite by our choice of $\ep_0$. Hence, we have proved (\ref{intconv}).\me

\no {\it 2.} Let $ \ep_0>0$ as in {\it 1.} We set $\displaystyle I_p=\int_{S_\si(\ep,T)\cap B} |\log|\De_\si(x_m\ga)|_\F^{-1}|^p d\ga.$\\

If $p$ is a positive integer, then there is positive constant $C'$ such that $\vert\log y\vert^p\leq C' y^{\ep_0/2}$ for all $y\geq 1$. 
Since $|\De_\si(x_m\ga)|^{-1}_\F\geq e^{\ep\Vert T\Vert}\geq 1$ for all  $\ga\in S_\si(\ep,T)$, we obtain
$$I_p\leq C'\int_{S_\si(\ep,T)\cap B} |\De_\si(x_m\ga)|^{-\ep_0/2}_\F d\ga\leq C' e^{-\frac{\ep\ep_0\Vert T\Vert}{2}}\int_B |\De_\si(x_m\ga)|^{-\ep_0}_\F d\ga.$$
If $p=0$ then by definition of  $  S_\si(\ep,T)$, one has  
 $$I_0=\int_{S_\si(\ep,T)\cap B} |\De_\si(x_m\ga)|_\F^{-\ep_0}|\De_\si(x_m\ga)|_\F^{\ep_0}d\ga\leq e^{- \ep\ep_0\Vert T\Vert}\int_{  B} |\De_\si(x_m\ga)|_\F^{-\ep_0}d\ga.$$
 In the two cases,  the result follows from {\it 1.}
\qed
\begin{lem}\label{KT+JT} Let $\ep_0>0$ as in Lemma \ref{intDe}. Given $\ep>0$, we can choose a constant $c>0$ such that for any $T\in a_{0,\F}$, one has 
$$\int_{S_\si(\ep,T)} \big( |K^T(x_m,\ga,f)|+|J^T(x_m,\ga,f)| \big)d\ga\leq ce^{-\frac{\ep\ep_0\Vert T\Vert}{4}}.$$

\end{lem}
\dem 
We recall that $$K^T(x_m, \ga,f)=|\De_\si(x_m\ga)|^{1/2}\int_{diag(A_M)\bb H\times H} \int_{diag(A_M)\bb H\times H} f_1(y_1^{-1}x_m\ga y_2) $$
$$\times f_2(x_1^{-1} x_m\ga   x_2) u_M(x_1,y_1,x_2,y_2, T)   d\overline{(x_1,x_2)} d\overline{(y_1,y_2)}$$ where 
$$u_M(x_1,y_1,x_2,y_2, T) =\int_{A_H\bb A_M} u(y_1^{-1}ax_1,T) u(y_2^{-1}ax_2, T) da.$$
We first establish an estimate of  $u_M$. Let $x,y\in H$ and $a\in A_M$. According to (\ref{Cartan}) applied to $H$, we can  write $y^{-1} ax= k_1 a_0k_2$ with $k_1, k_2\in K$ and $a_0\in A_0$. By definition of the norm, there is a positive constant $C_0$ such that 
$$ \log \Vert y^{-1} ax\Vert\leq C_0 (\Vert h_{A_0}(a_0)\Vert+1).$$

If  $u(y^{-1} ax)\neq 0$, then, by  definition of  $u(\cdot,T)$ (cf.  (\ref{uxT})),  the projection of $h_{A_0}(a_0)$ in $a_H\bb a_M$ belongs to the convex hull in $a_H\bb a_M$ of the $W(H,A_0)$-translates of $T$. Thus, there is a constant $C_1>0$ such that
\beq\label{inf1} \inf_{z\in A_H} \log \Vert y^{-1} za x\Vert\leq  C_1 (\Vert T\Vert+1).\eeq\\

We  assume that $\Vert T\Vert\geq 1$.  Taking $C_2=\max( 2C_1,1)$ and using  the property (\ref{prod})  of the norm,  we obtain

\beq\label{inflog} \inf_{z\in A_H} \log \Vert za \Vert\leq  C_2 (\Vert T\Vert+  \log\Vert x\Vert+\log  \Vert y\Vert).\eeq

We apply this  to $(x_1, y_1)$ and $(x_2, y_2)$ such that $u(y_1^{-1} ax_1,T) u(y_2^{-1} a x_2,T)\neq 0$. Hence, we deduce  that 
$$\inf_{z\in A_H} \log \Vert za\Vert \leq C_2 (\Vert T\Vert+\log\Vert x_1\Vert+\log  \Vert y_1\Vert+\log\Vert x_2\Vert+\log  \Vert y_2\Vert).$$

As $ \Vert x\Vert\leq \Vert x_m\Vert \Vert x_m^{-1} x\Vert$ and $1\leq \Vert T\Vert$, taking the integral over $a\in A_H\bb A_M$, we deduce the   following inequality
\ber\label{majuM}$$u_M(x_1,y_1,x_2,y_2, T) \preccurlyeq (\Vert T\Vert+\log\Vert  x_m^{-1}x_1\Vert+\log  \Vert  x_m^{-1}y_1\Vert+\log\Vert x_2\Vert+\log  \Vert y_2\Vert),$$
for all $x_1, y_1, x_2$ and $y_2$ in $H$. \eer
The function $u_M(x_1,y_1,x_2,y_2, T) $ is invariant by the diagonal (left) action of $A_M$ on $(x_1,x_2)$ and $(y_1,y_2)$ respectively. Since $x_m$ commutes with $A_S=A_M$ (cf.  Lemma \ref{fixeAS}), we can replace $\log\Vert x_m^{-1}x_1\Vert+ \log \Vert x_2\Vert$ and $\log\Vert  x_m^{-1}y_1\Vert+\log \Vert y_2\Vert$ by $\displaystyle \inf_{a\in A_M}\log\Vert (ax_m^{-1}x_1, ax_2)\Vert$ and $\displaystyle \inf_{a\in A_M}\log\Vert (ax_m^{-1}y_1, ay_2)\Vert$ respectively. By assumption, the quotient $A_M\bb S$ is compact, then, using (\ref{compact}),   one has
$$\inf_{a\in A_M}\Vert (a x_m^{-1}x, ax')\Vert\approx \inf_{s\in S}\Vert (sx_m^{-1}x, sx')\Vert,\quad x,x'\in H.$$\\
Therefore,  as $\Vert T\Vert \geq 1$, the inequality (\ref{majuM}) gives 
$$u_M(x_1,y_1,x_2,y_2, T) \preccurlyeq\Vert T\Vert+ \log \inf_{s\in S}\Vert (sx_m^{-1}x_1, sx_2)\Vert+ \log \inf_{s\in S}\Vert (sx_m^{-1}y_1, sy_2)\Vert, x_1, y_1, x_2, y_2\in H.$$
 In other words, this means that there are a positive constant $C_3$ and a positive integer $d$ such that, for all  $x_1, y_1, x_2$ and $y_2\in H$, one has 
$$u_M(x_1,y_1,x_2,y_2, T) \leq C_3(\Vert T\Vert+ \log \inf_{s\in S}\Vert (sx_m^{-1}x_1, sx_2)\Vert+ \log \inf_{s\in S}\Vert (sx_m^{-1}y_1, sy_2)\Vert)^d.$$

 Let $\Om$ be a compact set containing the support of $f_1$ and $f_2$. By   Lemma \ref{inf-delta}, there is a positive integer $k$ (independent of $\Om$) and a positive constant $C_\Om$ such that, if $x_m\ga\in x_m S_\si$ is a $\si$-regular point with $f_1(y_1^{-1}x_m\ga y_2) f_2(x_1^{-1} x_m\ga   x_2)\neq 0$ for some $x_1, x_2, y_1$ and $y_2$  in $H$  then 
$$u_M(x_1,y_1,x_2,y_2, T) \leq C_\Om(\Vert T\Vert+ \log\vert \De_\si(x_m\ga)\vert^{-k})^d.$$
This inequality  and the expression of $K^T(x_m, \ga,f)$ give
\beq\label{KTOI}\vert K^T(x_m, \ga,f)\vert \leq C_\Om(\Vert T\Vert+ \log\vert \De_\si(x_m\ga)\vert^{-k})^d\vert \cM(f_1)(x_m\ga)\cM(f_2)(x_m\ga)\vert,\eeq
where $\cM(f_j)$ is the  orbital integral of $f_j$ defined in (\ref{OI}). By Theorem \ref{OIbounded},  these orbital integrals are bounded by a constant $C_4$ on $ (x_mS_\si)\cap G^{\si-reg}.$ Hence, we obtain
$$\vert K^T(x_m, \ga,f)\vert \leq C_\Om C_4^2(\Vert T\Vert+ \log\vert \De_\si(x_m\ga)\vert^{-k})^d.$$

Let  $B$ be the set of $\ga$ in $S_\si$ such that $K^T(x_m,\ga,f)\neq 0$. Then $B$  is bounded by Theorem \ref{OIbounded} and (\ref{KTOI}). Using Lemma \ref{intDe},   we can find  a constant $C>0$ such that
\beq\label{majorKT} \int_{S_\si(\ep,T)}\vert K^T(x_m,\ga,f)\vert d\ga\leq C e^{-\frac{\ep\ep_0 \Vert T\Vert }{4}}.\eeq

If $\Vert T\Vert\leq 1$, then  (\ref{inf1}) implies that  if $u(x^{-1}a y)\neq 0$ then $$\inf_{z\in A_H} \log \Vert za \Vert\leq  2C_1 + \log\Vert  x\Vert+\log  \Vert y\Vert.$$
The same arguments to obtain (\ref{majuM})  imply   that   there is a positive constant $C'_1$ such that 
\beq\label{majuM1}u_M(x_1,y_1,x_2,y_2, T) \preccurlyeq (C'_1+\log\Vert  x_m^{-1}x_1\Vert+\log  \Vert  x_m^{-1}y_1\Vert+\log\Vert x_2\Vert+\log  \Vert y_2\Vert),\eeq
for $x_1, y_1, x_2$ and $y_2$ in $H$. Replacing $\Vert T\Vert $ by $C'_1$ in the reasoning after (\ref{majuM}), we deduce   that $\displaystyle \int_{S_\si(\ep,T)}\vert K^T(x_m,\ga,f)\vert d\ga $ is bounded. Hence, one obtains  (\ref{majorKT}) for $\Vert T\Vert\leq 1$.\me

We will now establish a similar estimate when $K^T$ is replaced by $J^T$. For this, it is enough to prove that the weight function $v_M$ have an estimate like (\ref{majuM}). We will see that this follows easily from the definition of $v_M$. Indeed,  for $x_1, y_1, x_2$ and $y_2$ in $H$, one has by definition
$$ v_M(x_1,y_1,x_2,y_2, T):= \int_{A_H\bb A_M}\si_M(h_M(a), \cY_M(x_1, y_1, x_2, y_2, T)) da$$
where $\si_M(\cdot, \cY_M(x_1, y_1, x_2, y_2, T))$ is a bounded function which vanishes in the complement of the convex hull $\cS_M(\cY_M(x_1, y_1, x_2, y_2, T))$ of the $(H,M)$-orthogonal set $ \cY_M(x_1, y_1, x_2, y_2, T)$  (cf. (\ref{SigmaMnul})). Since $ \cY_M(x_1, y_1, x_2, y_2, T)$ is the set of points $Z_P =\mbox{inf}^P(T_P+h_P(y_1)-h_{\overline{P}}(x_1),T_P+h_P(y_2)-h_{\overline{P}}(x_2))$ for $P\in\cP(M)$ (cf.  (\ref{YxyT})), if $\si_M(X, \cY_M(x_1, y_1, x_2, y_2, T))\neq 0$ then $\Vert X\Vert\leq \Vert Z_P\Vert$ for $P\in \cP(M)$. By definition of $T_P$, one has $\Vert T_P\Vert \leq \Vert T\Vert$.

Let us prove that for $P\in\cP(M)$, one has
\beq\label{majhP}\Vert h_P(x)\Vert\preccurlyeq 1+\log \Vert x\Vert, \quad x\in H.\eeq
We first  compare $\Vert m \Vert$ and $\Vert h_{M}(m)\Vert$  for $m\in M$. Let $M=K_MA_0K_M$ be the Cartan decomposition of $M$ where $K_M$ is a  suitable  compact subgroup of $M$. Then, each $m\in M$ can be written $m=ka(m)k'$ with $k,k'\in K_M$ and $a(m)\in A_0$.  Since $K_M$ is compact, the property   (\ref{compact}) gives   $\Vert m\Vert\approx \Vert a(m)\Vert,\ m\in M$ and this property does not depend on our choice of $a(m)$.  By (\ref{Normem0}), we have $\Vert a\Vert\approx e^{\Vert h_{A_0}(a)\Vert}, a\in A_0$. \\
Hence, there are a positive constant $C$ and a nonnegative integer $d$ such that $e^{\Vert h_{A_0}(a(m))\Vert}\leq C \Vert m\Vert^d$ for all $m\in M$. By (\ref{hGM}) applied to $(M, A_0)$, if  $a\in A_0$ then $h_{M}(a)$ is the orthogonal projection of $h_{A_0}(a)$ onto $a_{M}$, thus $\Vert h_{M}(a)\Vert\leq \Vert h_{A_0}(a)\Vert$. Since $h_{M}(m)=h_{M}(a(m))$ for $m\in M$, we obtain that  there is a  positive constant $C'$ such that
\beq\label{majhPM}\Vert h_{M}(m)\Vert\leq \Vert h_{A_0}(a(m))\Vert\leq C'(1+\log\Vert m\Vert),\quad m\in M.\eeq
By definition  (cf.  (\ref{mP}), (\ref{hP})), we have $h_P(x)=h_{M}(m_P(x)) $ for  $x\in H$   and by  (\ref{mn}), we have $\Vert m_P(x)\Vert\preccurlyeq \Vert x\Vert,  x\in H$. Thus, our claim (\ref{majhP}) follows from  (\ref{majhPM}).\me

Therefore, there are a positive $C_1$ and a positive integer $d$ such that if $\si_M(h_M(a), \cY_M(x_1, y_1, x_2, y_2, T))\neq 0$, then
$$\Vert h_M(a)\Vert\leq\Vert Z_P\Vert \leq C_1 (\Vert T\Vert+\log\Vert x_1\Vert+\log\Vert y_1\Vert+\log\Vert x_2\Vert+\log\Vert y_2\Vert)^d.$$
Since $\Vert x\Vert\leq\Vert x_m\Vert \Vert x_m^{-1} x\Vert  $ for $x\in H$, this gives the following estimates of $v_M$ analogous to (\ref{majuM}) and (\ref{majuM1}):
\ber\label{majvM}  If $\Vert T\Vert>1$ then

 $v_M(x_1,y_1,x_2,y_2, T) \preccurlyeq \Vert T\Vert+\log\Vert x_m^{-1}x_1\Vert+\log  \Vert x_m^{-1}y_1\Vert+\log\Vert x_2\Vert+\log  \Vert y_2\Vert,$ $x_1,y_1,x_2,y_2\in H$,\eer\\
and 
\ber\label{majvM1}there is a positive constant $C'_2$ such that for $\Vert T\Vert \leq 1$, one has 

$v_M(x_1,y_1,x_2,y_2, T) \preccurlyeq C'_2+\log\Vert x_m^{-1}x_1\Vert+\log  \Vert x_m^{-1}y_1\Vert+\log\Vert x_2\Vert+\log  \Vert y_2\Vert,$ $x_1,y_1,x_2,y_2\in H$.\eer
Arguing exactly as above for $K^T$, we deduce that there is a positive constant $C'$ such that 

$$\int_{S_\si(\ep,T)}\vert J^T(x_m,\ga,f)\vert d\ga\leq C' e^{-\frac{\ep\ep_0 \Vert T\Vert}{4}}.$$
This finishes the proof of the Lemma.
\qed

\begin{lem}\label{uM-vM} Fix $\de>0$. Then, there exist positive numbers $C,\ep_1$ and $\ep_2$ such that, for all $T$ with $d(T)\geq \de \Vert T\Vert$, and for all $x_1, y_1,x_2$ and $y_2$ in the set $H_{\ep_2}:= \{x\in H; \Vert x\Vert\leq e^{\ep_2\Vert T\Vert }\}$, one has
\beq\label{uM-vMestim}|u_M(x_1,y_1,x_2,y_2, T)-v_M(x_1,y_1,x_2,y_2,T)|\leq C e^{-\ep_1 \Vert T\Vert}.\eeq
\end{lem}

\dem  If $\Vert T\Vert $ remains bounded then, by (\ref{majuM}), (\ref{majuM1}),  (\ref{majvM}) and (\ref{majvM1}),  the functions $u_M$ and $v_M$ are  bounded and  the result (\ref{uM-vMestim}) is trivial. Thus we have to prove the Lemma for  $\Vert T\Vert$  sufficiently large  and $d(T)\geq \de \Vert T\Vert$. 

 By  (5.8) of \cite{ArLT}, we can choose $\ep_2$ such that   $d(\cY_M(x,y,T))>0$   for all $x,y\in H_{\ep_2}$. By the discussion of   l.c. bottom of page 38 and top of page 39,  there is a constant $C_0>0$ such that, for $T$ with $d(T)\geq \de \Vert T\Vert$ and $\Vert T\Vert> C_0$, for  $x,y\in H_{\ep_2}$ and $a\in A_H\bb A_M$, one has 
$$u(y^{-1} a x,T)= \si_M( h_M(a), \cY_M(x,y,T)).$$
By Lemma \ref{infY}, for $X\in a_M$, we have $$\si_M( X, \cY_M(x_1,y_1, x_2, y_2,T))=\si_M( X, \cY_M(x_1,y_1,T))\si_M( X, \cY_M(x_2,y_2,T)).$$\\
Thus, one deduces that 
$$\si_M( h_M(a), \cY_M(x_1,y_1, x_2, y_2,T))=u(y_1^{-1} a x_1,T)u(y_2^{-1} a x_2,T),$$
 for $a\in A_H\bb A_M$.  Hence, for $d(T)\geq \de \Vert T\Vert\geq \de C_0$, and $x_i, y_i$ in $H_{\ep_2}$, we have 
$$u_M(x_1,y_1,x_2,y_2, T)=v_M(x_1,y_1,x_2,y_2,T).$$
This finishes the proof of  the Lemma. \qed

\no Theorem \ref{resultatprincipal} follows from the   corollary below.
\begin{cor}\label{KTx-JTx}  Fix $\de>0$.There exist two positive numbers $\ep$ and $c>0$ such that,  for all $T$ with $d(T)\geq \de \Vert T\Vert$, one has

\beq\int_{\ga\in S_\si}   |K^T(x_m, \ga,f) -J^T(x_m\ga,f)| \ d\ga\leq c e^{- \ep\Vert T\Vert}.\eeq

\end{cor} 
\dem By Lemma \ref{KT+JT}, it is enough to prove that  we can find positive numbers  $\ep$, $\ep'$ and $C_0$ such that
\beq\label{estimKT-JT}\int_{\ga\in S_\si-S_\si(\ep,T)}   |K^T(x_m, \ga,f)-J^T(x_m,\ga,f)|  d\ga\leq C_0e^{-  \ep'\Vert T\Vert}\eeq

\no where $S_\si(\ep,T)$ is defined in (\ref{SsepT}).

Let $\ep>0$.
Let $\Om$ be a compact subset of $G$ which contains  the support of $f_1$ and $f_2$. We will estimate   $\vert u_M(x_1,y_1, x_2,y_2,T)-v_M(x_1,y_1, x_2,y_2,T)\vert$ for $x_1,x_2, y_1$ and $y_2$ in $H$ satisfying $ x_1^{-1} x_m \ga x_2\in \Om$ and $y_1^{-1} x_m\ga  y_2\in\Om$ for some $\ga\in S_\si-S_\si(\ep,T)$ with $x_m\ga\in G^{\si-reg}$. For this, we will use the invariance of the functions $u_M$ and $v_M$   by the diagonal left action of $A_M$ on $(x_1,x_2)$ and  $(y_1,y_2)$ respectively. 

 By Lemma \ref{inf-delta}, there are  a positive integer $k$ and a positive constant  $C_\Om$, (depending only on $\Om$) such that, for all $\ga\in S_\si-S_\si(\ep, T)$ with $x_m\ga\in G^{\si-reg}$ and  for all $x_i,y_i$ in $H$, $i=1,2$  with $ x_1^{-1} x_m \ga x_2$ and $y_1^{-1} x_m\ga  y_2$ in $\Om$, we have 
\beq\label{infx}\inf_{s\in S}\Vert (sx_m^{-1}x_1, sx_2)\Vert\leq C_\Om \De_\si(x_m\ga)^{-k}\leq C_\Om e^{k\ep \Vert T\Vert}\eeq
and 
$$\inf_{s\in S}\Vert (sx_m^{-1}y_1, sy_2)\Vert\leq C_\Om \De_\si(x_m\ga)^{-k}\leq C_\Om e^{k\ep \Vert T\Vert}.$$\\
 Since $A_M\bb S$ is compact, we deduce from (\ref{compact}) and (\ref{infx}) that there is a constant $C'_\Om>0$ such that
 $$\inf_{a\in A_M}\Vert (ax_m^{-1}x_1, ax_2)\Vert\leq C'_\Om e^{k\ep \Vert T\Vert}.$$\\
Thus,  for $\eta>0$, there is $a_0\in A_M$ such that 
\beq\label{majx}\Vert a_0 x_m^{-1} x_1\Vert \Vert a_0 x_2\Vert\leq C_\Om e^{k\ep\Vert T\Vert}+\eta.\eeq
Since $A_M=A_S$, the point $a_0$ commutes with $x_m$ by (\ref{xm}) and   we have $\Vert a_0x_1\Vert\leq \Vert x_m\Vert \Vert x_m^{-1}  a_0 x_1\Vert$.\me

If $\Vert T\Vert$ remains bounded, then $\Vert a_0x_i\Vert, i=1,2$ are bounded by a constant independent of $\Vert T\Vert$. By the same arguments, there is $a_1\in A_M$ such that $\Vert a_1y_i\Vert, i=1,2$ are bounded by a constant independent of $\Vert T\Vert$. Using the invariance of $u_M$ and $v_M$ by the left action of $diag(A_M)$ on $(x_1,x_2)$ and $(y_1,y_2)$ respectively and the estimates (\ref{majuM}), (\ref{majuM1}) , (\ref{majvM}) and (\ref{majvM1}) for $u_M$ and $v_M$, we deduce that $\vert u_M(x_1,y_1, x_2,y_2,T)-v_M(x_1,y_1, x_2,y_2,T)\vert$ is bounded by a constant independent of $T$ and of $x_i, y_i$.
Recall that by Theorem \ref{OIbounded}, the constant
$$C_1:= \int_{S_\si} \cM(\vert f_1\vert)(x_m\ga)\cM(\vert f_2\vert)(x_m\ga) d\ga$$
is finite. We deduce that $\int_{\ga\in S_\si-S_\si(\ep,T)}   |K^T(x_m, \ga,f)-J^T(x_m,\ga,f)|  d\ga$ is bounded, hence we obtain (\ref{estimKT-JT}).\me

We assume that $\Vert T\Vert$ is not bounded. Let $\ep_1,\ep_2$ and $C$ as in    Lemma \ref{uM-vM}.  Taking $\Vert T\Vert $ to be sufficiently large and $\ep$ such that $k\ep$ is smaller than the constant $\ep_2$, we can assume by (\ref{majx}) that 
$$\Vert a_0x_i\Vert\leq e^{\ep_2\Vert T\Vert}, \quad i=1,2.$$\\
The same arguments are valid for  $\Vert y_i\Vert $, $i=1,2$. Thus, there is $a_1\in A_M$ such that 
  $$\Vert a_1y_i\Vert\leq e^{\ep_2\Vert T\Vert}, \quad i=1,2.$$\\
Using  Lemma \ref{uM-vM} and the invariance  of $u_M$ and $v_M$ by the left action of the diagonal of $A_M$ on $(x_1, x_2)$ and $(y_1, y_2)$ respectively, we deduce that, for all $T$ with $d(T)\geq \de \Vert T\Vert$, one has
$$\vert u_M(x_1,y_1, x_2,y_2,T)-v_M(x_1,y_1, x_2,y_2,T)\vert \leq C e^{-\ep_1\Vert T\Vert}.$$\\
Hence, we obtain
$$\int_{S-S_\si(\ep,T)}\vert K^T(x_m,\ga,f)-J^T(x_m,\ga,T)\vert\leq C C_1 e^{-\ep_1\Vert T\vert},$$
where $C_1:= \int_{S_\si} \cM(\vert f_1\vert)(x_m\ga)\cM(\vert f_2\vert)(x_m\ga) d\ga$.
This finishes the proof of the Corollary.\qed


\subsection{ The function $J^T(f)$}
The goal of this section is to prove that  $J^T(f)$ is of the form
\beq\label{expressJT}\sum_{k=0}^N p_k(T,f) e^{\xi_k(T)},\eeq
where  $\xi_0=0, \xi_1,\ldots, \xi_N$ are distinct points in $i a_0^*$ and each $p_k(T,f)$ is a polynomial function of $T$. Moreover, the constant term $\tilde{J}(f):=p_0(0,f)$ is well-defined and is uniquely determined by $K^T(f)$. Except for one detail, our arguments and calculations are the same as  those of section 6 of \cite{ArLT}. We give the details of proof for convenience of the reader.\me

Recall that $J^T(f)$ is a finite sum of the distributions 
$$J^T(x_m, \ga,f)=|\De_\si(x_m\ga)|_\F^{1/2}\int_{diag(A_M)\bb H\times H} \int_{diag(A_M)\bb H\times H} f_1(y_1^{-1}x_m\ga y_2) $$
$$\times f_2(x_1^{-1} x_m\ga   x_2) v_M(x_1,y_1,x_2,y_2, T)   d\overline{(x_1,x_2)} d\overline{(y_1,y_2)}$$
where $M \in \cL(A_0)$, $S$ is a maximal torus of $M$ such that $A_S=A_M$, $x_m\in\kappa_{S}$ and $v_M(x_1,y_1,x_2,y_2, T) := \int_{A_H\bb A_M}\si_M(h_M(a), \cY_M(x_1, y_1, x_2, y_2, T)) da$ where $\cY_M(x_1, y_1, x_2, y_2, T)$ is  defined in (\ref{YxyT}).\me

We first study the weight function $v_M$ as a function of $T$. We fix $M\in\cL(A_0)$ and  $x_1,y_1, x_2$ and $y_2$ in $H$. 

Let 
$\Lr_M:=(a_{M,\F}+a_H)/a_H$ and $
\widetilde{\Lr}_M:=(\tilde{a}_{M,\F}+a_H)/ a_{H}$ be the projection in $a_M/a_H$ of the lattices $a_{M,\F}$ and $\tilde{a}_{M,\F}$ respectively.
By (\ref{projreseau}), one has 
\beq\label{Ltilde}\tilde{a}_{M,\F}/ \tilde{a}_{H,\F}=\tilde{a}_{M,\F} /\tilde{a}_{M,\F}\cap a_H\simeq \widetilde{\Lr}_M .\eeq\\
For $M=A_0$, we replace the subscript $A_0$ by $0$. We denote by  $\Lr^\vee:=\mbox{Hom}(\Lr, 2\pi i\Z)$ the dual lattice  of  a lattice $\Lr$. 

Let $P\in\cP(M)$. We introduce the following sublattice of $\Lr_M$. For $k\in\N$, we set 
$$\mu_{\al,k}:=k\log (q) \check{\al}, \al\in \De_P,$$
where $q$ is the order of the residual field of $\F$, 
and 
$$\Lr_{M,k}:=\sum_{\al\in \De_P} \Z\mu_{\al,k}.$$
Then $\Lr_{M,k}$ is a lattice  in $a_M^H\simeq a_M/a_H$ independent of $P$  and by (\cite{ArWOI} \textsection 4), one can find $k\in\N^*$ such that for all $M\in \cL(A_0)$, one has
$$\Lr_{M,k}\subset \widetilde{\Lr}_M.$$\\
The set of points  $\sum_{\al\in\De_P} y_\al \mu_{\al,k}$ with $y_\al\in]-1,0] $ is  a fundamental domain of $\Lr_{M,k}$ which we denote by $\cD_{M,k}$. 
\ber\label{XPY} For $X\in\Lr_M/\Lr_{M,k}$ and $Y\in  a_M/a_H$, we denote by $\bar{X}_P(Y)$ the representative of $X$ in $\Lr_M$ such that 
$\bar{X}_P(Y)-Y\in \cD_{M,k}$.\eer\\
For $\la\in a_{M, \C}^*$, we set 
\beq\label{theta} \te_{P,k}(\la)=vol(a_M^H/\Lr_{M,k})^{-1}\prod_{\al\in \De_P}(1-e^{-\la(\mu_{\al, k})}).\eeq

We fix  $T\in a_{0,\F}$. By definition of $\si_M$  (cf.  (\ref{SigmaM})), the function $v_M$ depends only on the image of $T_P$ in $\Lr_M$.  Hence we can assume that $T$ lies in the lattice $\Lr_0$. For $P\in\cP(M)$, the map $T\mapsto T_P$ sends surjectively $\Lr_0$ onto the intersection of $\Lr_M$ with the closure $\overline{a_P^+}$ of the chamber associated to $P$. Thus, we may restrict $T$ to lie in the intersection of $\Lr_0$ with suitable regular points in some positive chamber $a_0^+$ of $a_H\bb a_0$. Then  the points $T_P$ range over a suitable regular points in $\Lr_M\cap a_P^+$. \me


We recall that $ \cY_M(x_1, y_1, x_2, y_2, T)$ is the set of points $Z_P:=Z_P(x_1, y_1, x_2, y_2, T)$ defined in (\ref{ZP}).  Thus, we can write 
\beq\label{ZP} Z_P=T_P+Z_P^0\mbox{ with }Z^0_P:=\mbox{inf}^P(h_P(y_1)-h_{\overline P}(x_1), h_P(y_2)-h_{\overline P}(x_2)).\eeq

Notice that the points $Z_P^0$ do not necessarily belong to the lattice $\Lr_M$. It is the only difference with \cite{ArLT} section 6 in what follows.
  \begin{lem}\label{dvptvM} There is a positive integer $N$ independent of $M$ and polynomial functions $q_\xi(T)$ for $\xi\in \bi(\frac{1}{N}\Lr_0^\vee\big)/\Lr_0^\vee$ (depending on $x_1, y_1, x_2$ and $y_2$), such that
 
 $$v_M(x_1, y_1, x_2, y_2, T)=\sum_{\xi\in (\frac{1}{N}\Lr_0^\vee)/\Lr_0^\vee} q_\xi(T)e^{\xi(T)}.$$
 Moreover, the constant term $\tilde{v}_M(x_1,y_1,x_2,y_2):=q_0(0)$ of $v_M(x_1, y_1, x_2, y_2, T)$ is given by
 $$\tilde{v}_M(x_1,y_1,x_2,y_2)=\lim_{\La\to0}\big(\sum_{P\in\cP(M)}\vert \Lr_M/\Lr_{M,k}\vert^{-1}\sum_{X\in\Lr_M/\Lr_{M,k}}e^{\langle \La,\bar{X}_P(Z_P^0)\rangle} \te_{P,k}(\La)^{-1}\big).$$
 \end{lem}
 \dem 
The kernel of the surjective map $h_M: A_H\bb A_M\to \tilde{a}_{M,\F}/ \tilde{a}_{H,\F}$ is a compact group which has volume $1$ by our convention of choice of measure. Thus, using (\ref{Ltilde}), we can write
$$v_M(x_1,y_1,x_2,y_2, T) :=\sum_{X\in \widetilde{\Lr}_M}\si_M(X, \cY_M(x_1, y_1, x_2, y_2, T)).$$\\
For our study, it is convenient to take a  sum over $\Lr_M$. The finite quotient $\widetilde{\Lr}_M^\vee/\Lr_M^\vee$ can be identified with the character group of $\Lr_M/\widetilde{\Lr}_M$ under the pairing $$(\nu,X)\in \widetilde{\Lr}_M^\vee/\Lr_M^\vee\times \Lr_M/\widetilde{\Lr}_M\mapsto e^{\nu(X)}.$$
 Hence, by inversion formula on finite abelian groups,  we obtain
$$v_M(x_1,y_1,x_2,y_2, T)=\vert \Lr_M/\widetilde{\Lr}_M\vert^{-1}\sum_{\nu\in  \widetilde{\Lr}_M^\vee/\Lr_M^\vee}\sum_{X\in  \Lr_M} \si_M(X,  \cY_M(x_1, y_1, x_2, y_2, T)) e^{\nu(X)}.$$\\
Coming back to the  definition of $\si_M$ (cf.  (\ref{SigmaM})), we fix a small point $\La\in (a_M/a_H)^*_\C$ whose real part $\La_R$ is in general position. One has
$$ \si_M(X,  \cY_M(x_1, y_1, x_2, y_2, T))=\sum_{P\in\cP(M)} (-1)^{\vert\De_P^\La\vert}\varphi_P^\La (X-Z_P)$$
$$= \lim_{\La\to 0}\sum_{P\in\cP(M)} (-1)^{\vert\De_P^\La\vert}\varphi_P^\La (X-Z_P)e^{\La(X)}.$$\\
By definition of $\varphi_P^\La$, the function $X\mapsto e^{\La (X)}$ is rapidly decreasing on the support of $X\mapsto \varphi_P^\La(X-Z_P)$. Hence, the product of these two functions is summable over $X\in\Lr_M$. Therefore, we can write 
\beq\label{vMFP}v_M(x_1,y_1,x_2,y_2, T) =\sum_{\nu\in  \widetilde{\Lr}_M^\vee/\Lr_M^\vee}\lim_{\La\to 0}
\sum_{P\in\cP(M)} F_P^T(\La)\eeq
where
$$F_P^T(\La):= \vert \Lr_M/\widetilde{\Lr}_M\vert^{-1}\sum_{X\in  \Lr_M}(-1)^{\vert\De_P^\La\vert}\varphi_P^\La (X-Z_P)e^{(\La+\nu)(X)}.$$
The above discussion implies that 
\ber\label{FPana} the map $\La\mapsto \sum_{P\in\cP(M)} F_P^T(\La)$ is analytic at $\La=0$.\eer 
We fix $P\in\cP(M)$. We want to express $F_P^T(\La)$  in terms of a product of geometric series. For this, we write
\beq\label{FP1}F_P^T(\La):= \vert \Lr_M/\widetilde{\Lr}_M\vert^{-1}\sum_{X\in  \Lr_M/\Lr_{M,k}}\sum_{X'\in\Lr_{M,k}}(-1)^{\vert\De_P^\La\vert}\varphi_P^\La (X+X'-Z_P)e^{(\La+\nu)(X+X')}.\eeq\\
Let $X\in\Lr_M/\Lr_{M,k}$. Recall that $\bar{X}_P(Y) $ is the representative of $X$ in $\Lr_M$ such that $\bar{X}_P(Y)-Y\in\cD_{M,k}$. We set 
$$\bar{X}_P^\La(Y):=\bar{X}_P(Y)+\sum_{\al\in\De_P^\La}\mu_{\al,k}.$$\\
Thus $\bar{X}_P^\La(Y)$ is also a representative of $X$ in $\Lr_M$. Taking $Y:=Z_P$, we can set 
$$\varphi_P^\La (X +X'-Z_P)=\varphi_P^\La (\bar{X}_P^\La(Z_P)+X'-Z_P) $$
in (\ref{FP1}). The set of points $X'\in\Lr_{M,k}$ such that this characteristic function equals to $1$ is exactly the set 
$$\{\sum_{\al\in \De_P^\La} n_\al\mu_{\al,k}-\sum_{\al\in \De_P-\De_P^\La} n_\al\mu_{\al,k}; n_\al\in\N\}.$$ \\
Therefore, a simple  calculation as in \cite{ArLT} top of page 45 gives
\ber\label{egalite}$$(-1)^{\vert\De_P^\La\vert}\sum_{X'\in\Lr_{M,k}}\varphi_P^\La (X+X'-Z_P)e^{(\La+\nu)(X+X')}$$
$$= e^{(\La+\nu)(\bar{X}_P(Z_P))}(\prod_{\al\in \De_P}(1-e^{-(\La+\nu)(\mu_{\al,k})})^{-1}.$$\eer\\
We have fixed the Haar measure on $a_M^H\simeq a_M/a_G$ with the property that the quotient of $a_M/a_H$ by the lattice $\widetilde{\Lr}_M$ has volume $1$. Thus, we have
 $$\vert \Lr_M/\widetilde{\Lr}_M\vert^{-1} \prod_{\al\in \De_P}(1-e^{-(\La+\nu)(\mu_{\al,k})})^{-1} =\vert \Lr_M/\widetilde{\Lr}_{M,k}\vert^{-1}\te_{P,k}(\La+\nu)^{-1}.$$
 
By the above equality, (\ref{FP1}) and (\ref{egalite}), we obtain 
\beq\label{FP2}F_P^T(\La)=\vert\Lr_M/\Lr_{M,k}\vert^{-1}\sum_{X\in \Lr_M/\Lr_{M,k}} e^{<\La+\nu, \bar{X}_P(Z_P)>}\te_{P,k}(\La+\nu)^{-1}.\eeq

Let $X\in \Lr_M/\Lr_{M,k}$. We recall that $T_P$ belongs to $\Lr_M$ for $P\in\cP(M)$ and $Z_P=T_P+Z_P^0$ (cf. (\ref{ZP})). By definition (cf.  (\ref{XPY})), the point $\bar{X}_P(Z_P)$ is the unique representative of $X$ in $\cL_M$ such that $ \bar{X}_P(Z_P)- T_P-Z_P^0\in \cD_{M,k}$ and $\overline{(X-T_P)}_P(Z^0_P)$ is the unique representative of $X-T_P$ in $\Lr_M$ such that $\overline{(X-T_P)}_P(Z^0_P)-Z^0_P\in\cD_{M,k}$. Hence, we deduce that
\beq\label{XZP} \bar{X}_P(Z_P)= \overline{(X-T_P)}_P(Z^0_P)+T_P.\eeq

Replacing $X$ by $X-T_P$ in (\ref{FP2}), we obtain 
\beq\label{FP3}F_P(\La)^T= \vert\Lr_M/\Lr_{M,k}\vert^{-1}\sum_{X\in \Lr_M/\Lr_{M,k}} e^{<\La+\nu, T_P+\bar{X}_P(Z^0_P)>}\te_{P,k}(\La+\nu)^{-1}\eeq
where $\bar{X}_P(Z^0_P)$ is independent of $T$. Thus by (\ref{vMFP}), we have established that 
$v_M(x_1,y_1,x_2,y_2,T)$ is equal to
\beq\label{vMfinal}\sum_{\nu\in  \widetilde{\Lr}_M^\vee/\Lr_M^\vee}\lim_{\La\to 0} \big(\sum_{P\in\cP(M)}\vert\Lr_M/\Lr_{M,k}\vert^{-1}\sum_{X\in \Lr_M/\Lr_{M,k}} e^{<\La+\nu, T_P+\bar{X}_P(Z^0_P)>}\te_{P,k}(\La+\nu)^{-1}\big).\eeq\\
Recall that the expression in the brackets is analytic at $\La=0$ (cf.  (\ref{FPana})). To analyze this expression as function of $T$, we argue as in (\cite{W1} p.315). We give the details for convenience of lecture. We replace $\La$ by $z\La$. The map $z\mapsto \te_{P,k}(z\La+\nu)^{-1}$ may have a pole at $z=0$. Let $r$ denotes the biggest order of this pole when $P$ runs $\cP(M)$. Then, using Taylor expansions, one deduces that 
$$\lim_{\La\to 0} \big(\sum_{P\in\cP(M)}\vert\Lr_M/\Lr_{M,k}\vert^{-1}\sum_{X\in \Lr_M/\Lr_{M,k}} e^{<\La+\nu, T_P+\bar{X}_P(Z^0_P)>}\te_{P,k}(\La+\nu)^{-1}\big)=$$
$$\sum_{m=0}^r\sum_{P\in\cP(M)}C_m\sum_{X\in \Lr_M/\Lr_{M,k}}\frac{\partial^m}{\partial z^m}(e^{<z\La+\nu, T_P+\bar{X}_P(Z^0_P)>})_{[z=0]} \frac{\partial^{r-m}}{\partial z^{r-m}}(z^r\te_{P,k}(z\La+\nu)^{-1})_{[z=0]} ,$$
where $C_m=\frac{1}{m! (r-m)!}\vert\Lr_M/\Lr_{M,k}\vert^{-1}$. 

But we have
$$\frac{\partial^m}{\partial z^m}(e^{<z\La+\nu, T_P+\bar{X}_P(Z^0_P)>})_{[z=0]} =(<\La,  T_P+\bar{X}_P(Z^0_P)>)^m e^{<\nu, T_P+\bar{X}_P(Z^0_P)>},$$
and $\dfrac{\partial^{r-m}}{\partial z^{r-m}}(z^r\te_{P,k}(z\La+\nu)^{-1})_{[z=0]}$ is independent of $T_P$.\\
Therefore, we deduce that $v_M(x_1,y_1,x_2,y_2, T)$ is a finite sum of functions
$$q_{P,\nu}(T_P) e^{\nu(T_P)}, \quad \nu\in \widetilde{\Lr}_M^\vee/\Lr_M^\vee, P\in\cP(M),$$
where $q_{P,\nu}$ is a polynomial function on $a_M$.\me

Since $ \Lr_0^\vee\subset  \widetilde{\Lr_0}^\vee$ are lattices of same rank, one can find a positive integer $N$ such that $N\widetilde{\Lr_0}^\vee\subset \Lr_0^\vee$. Therefore, by our choice of $T$ 
 and the above expression,   we can write $$v_M(x_1,y_1,x_2,y_2, T)=\sum_{\xi\in (\frac{1}{N}\Lr_0^\vee )/\Lr_0^\vee} q_\xi(T)e^{\xi(T)},$$
 where $q_\xi(T)$ is a polynomial function of $T$. 
 This gives the first part of the Lemma.
 
Since the polynomials $q_\xi(T)$ are obviously uniquely determined, the constant term $\tilde{v}_M(x_1,y_1,x_2,y_2):=q_0(0)$ is well defined. To calculate it, we take the summand corresponding to $\nu=0$ in (\ref{vMfinal}) 	and then set $T=0$. We obtain
$$\tilde{v}_M(x_1,y_1,x_2,y_2)=\lim_{\La\to0}\big(\sum_{P\in\cP(M)}\vert \Lr_M/\Lr_{M,k}\vert^{-1}\sum_{X\in\Lr_M/\Lr_{M,k}}e^{\langle \La,\bar{X}_P(Z_P^0)\rangle} \te_{P,k}(\La)^{-1}\big).$$
This finishes the proof of the Lemma.\qed

We substitute the expression we have obtained for  $v_M$  in Lemma \ref{dvptvM} into the expression for $J^T(x_m, \ga, f)$. Hence, we obtain the following  similar decomposition for $J^T(f)$.
 \begin{cor}\label{expressionJT} There is a decomposition
$$J^T(f)=\sum_{\xi\in (\frac{1}{N}\Lr_0^\vee )/\Lr_0^\vee} p_\xi(T,f)e^{\xi(T)},\quad T\in\Lr_0\cap a_0^+,$$
 where $N$ is positive integer and each $p_\xi(T,f)$ is a polynomial function of $T$. Moreover, the constant term $\Tilde{J}(f):=p_0(0,f)$ of $J^T(f)$ is given by
 $$\tilde{J}(f)= J^T(f):=\sum_{M\in\cL(A_0)} c_M\sum_{S\in\cT_M}\sum_{x_m\in\kappa_S} c_{S,x_m}\int_{S_\si} \tilde{J}(x_m,\ga,f)d\ga,$$
 where
 $$ \tilde{J}(x_m,\ga,f)=|\De_\si(x_m\ga)|^{1/2}\int_{diag(A_M)\bb H\times H} \int_{diag(A_M)\bb H\times H} f_1(y_1^{-1}x_m\ga y_2) $$
$$\times f_2(x_1^{-1} x_m\ga   x_2) \tilde{v}_M(x_1,y_1,x_2,y_2)   d\overline{(x_1,x_2)} d\overline{(y_1,y_2)}.$$
\end{cor}

\noindent P.~Delorme, Aix-Marseille Université, CNRS, Centrale Marseille, I2M, UMR 7373, 13453
Marseille, France.\\
{\it E-mail address}:  patrick.delorme@univ-amu.fr\me

\noindent P.~Harinck, Ecole polytechnique,  CMLS - CNRS, UMR 7640, 
Route de Saclay, 91128 Palaiseau Cédex, France.\\  {\it E-mail address}: pascale.harinck@ polytechnique.edu\me

\noindent S.~Souaifi,  Université de Strasbourg, IRMA CNRS, UMR 7501, 7 rue Descartes, 67084 Strasbourg Cédex, France.\\
{\it E-mail address}:  sofiane.souaifi@math.unistra.fr\\

\end{document}